# ON IMPROVEMENT IN ESTIMATING POPULATION PARAMETER(S) USING AUXILIARY INFORMATION

**Rajesh Singh**

**Department of Statistics, BHU, Varanasi (U.P.), India**

**Florentin Smarandache**

**Chair of Department of Mathematics, University of New Mexico, Gallup, USA**

**Preface**

The purpose of writing this book is to suggest some improved estimators using auxiliary information in sampling schemes like simple random sampling and systematic sampling.

This volume is a collection of five papers. The following problems have been discussed in the book:

In chapter one an estimator in systematic sampling using auxiliary information is studied in the presence of non-response. In second chapter some improved estimators are suggested using auxiliary information. In third chapter some improved ratio-type estimators are suggested and their properties are studied under second order of approximation.

In chapter four and five some estimators are proposed for estimating unknown population parameter(s) and their properties are studied.

This book will be helpful for the researchers and students who are working in the field of finite population estimation.

**Contents**

# Use of Auxiliary Information for Estimating Population Mean in Systematic Sampling under Non- Response

**[1]Manoj K. Chaudhary, [1]Sachin Malik, †[1]Rajesh Singh and [2]Florentin Smarandache**

**[1]Department of Statistics, Banaras Hindu University,Varanasi-221005, India**

**[2]Chair of Department of Mathematics, University of New Mexico, Gallup, USA**

***†Corresponding author***

**Abstract**

In this paper we have adapted Singh and Shukla (1987) estimator in systematic sampling using auxiliary information in the presence of non-response. The properties of the suggested family have been discussed. Expressions for the bias and mean square error (MSE) of the suggested family have been derived. The comparative study of the optimum estimator of the family with ratio, product, dual to ratio and sample mean estimators in systematic sampling under non-response has also been done. One numerical illustration is carried out to verify the theoretical results.

## 1. Introduction

There are some natural populations like forests etc., where it is not possible to apply easily the simple random sampling or other sampling schemes for estimating the population characteristics. In such situations, one can easily implement the method of systematic sampling for selecting a sample from the population. In this sampling scheme, only the first unit is selected at random, the rest being automatically selected according to a predetermined pattern. Systematic sampling has been considered in detail by Madow and Madow (1944), Cochran (1946) and Lahiri (1954). The application of systematic sampling to forest surveys has been illustrated by Hasel (1942), Finney (1948) and Nair and Bhargava (1951).

The use of auxiliary information has been permeated the important role to improve the efficiency of the estimators in systematic sampling. Kushwaha and Singh (1989) suggested a class of almost unbiased ratio and product type estimators for estimating the population mean using jack-knife technique initiated by Quenouille (1956). Later Banarasi et al. (1993), Singh and Singh (1998), Singh et al. (2012), Singh et al. (2012) and Singh and Solanki (2012) have made an attempt to improve the estimators of population mean using auxiliary information in systematic sampling.

The problem of non-response is very common in surveys and consequently the estimators may produce bias results. Hansen and Hurwitz (1946) considered the problem of estimation of population mean under non-response. They proposed a sampling plan that involves taking a subsample of non-respondents after the first mail attempt and then enumerating the subsample by personal interview. El-Badry (1956) extended Hansen and Hurwitz (1946) technique. Hansen and Hurwitz (1946) technique in simple random sampling is described as: From a population $U = (U_1, U_2, ---, U_N)$, a large first phase sample of size n' is selected by simple random sampling without replacement ( SRSWOR). A smaller second phase sample of size n is selected from n' by SRSWOR. Non-response occurs on the second phase of size n in which $n_1$ units respond and $n_2$ units do not. From the $n_2$ non-respondents, by SRSWOR a sample of $r = n_2/ k$; $k > 1$units is selected. It is assumed that all the r units respond this time round. ( see Singh and Kumar (20009)). Several authors such as Cochran (1977), Sodipo and Obisesan ( 2007), Rao (1987), Khare and Srivastava ( 1997) and Okafor and Lee (2000) have studied the problem of non-response under SRS.

In the sequence of improving the estimator, Singh and Shukla (1987) proposed a family of factor-type estimators for estimating the population mean in simple random sampling using an auxiliary variable, as

$$T_\alpha = \bar{y}\left[\frac{(A+C)\bar{X} + fB\bar{x}}{(A+fB)\bar{X} + C\bar{x}}\right] \tag{1.1}$$

where $\bar{y}$ and $\bar{x}$ are the sample means of the population means $\bar{Y}$ and $\bar{X}$ respectively. $A$, $B$ and $C$ are the functions of $\alpha$, which is a scalar and chosen so as the MSE of the estimator $T_\alpha$ is minimum.

Where,

$$A = (\alpha-1)(\alpha-2),\ B = (\alpha-1)(\alpha-4),$$

$$C = (\alpha-2)(\alpha-3)(\alpha-4);\ \alpha > 0 \text{ and}$$

$$f = \frac{n}{N}.$$

**Remark 1** : If we take $\alpha$ = 1, 2, 3 and 4, the resulting estimators will be ratio, product, dual to ratio and sample mean estimators of population mean in simple random sampling respectively (for details see Singh and Shukla (1987) ).

In this paper, we have proposed a family of factor-type estimators for estimating the population mean in systematic sampling in the presence of non-response adapting Singh and Shukla (1987) estimator. The properties of the proposed family have been discussed with the help of empirical study.

## 2. Sampling Strategy and Estimation Procedure

Let us assume that a population consists of $N$ units numbered from 1 to $N$ in some order. If $N = nk$, where $k$ is a positive integer, then there will be $k$ possible samples each consisting of $n$ units. We select a sample at random and collect the information from the units of the selected sample. Let $n_1$ units in the sample responded and $n_2$ units did not respond, so that $n_1 + n_2 = n$. The $n_1$ units may be regarded as a sample from the response class and $n_2$ units as a sample from the non-response class belonging to the population. Let us assume that $N_1$ and $N_2$ be the number of units in the response class and non-response class respectively in the population. Obviously, $N_1$ and $N_2$ are not known but their unbiased estimates can be obtained from the sample as

$$\hat{N}_1 = n_1 N / n; \;\; \hat{N}_2 = n_2 N / n.$$

Further, using Hansen and Hurwitz (1946) technique we select a sub-sample of size $h_2$ from the $n_2$ non-respondent units such that $n_2 = h_2 L$ ($L > 1$) and gather the information on all the units selected in the sub-sample (for details on Hansen and Hurwitz (1946) technique see Singh and Kumar (2009)).

Let $Y$ and $X$ be the study and auxiliary variables with respective population means $\overline{Y}$ and $\overline{X}$. Let $y_{ij}(x_{ij})$ be the observation on the $j^{th}$ unit in the $i^{th}$ systematic sample under study (auxiliary) variable ($i = 1...k : j = 1...n$).Let us consider the situation in which non-response is observed on study variable and auxiliary variable is free from non-response. The Hansen-Hurwitz (1946) estimator of population mean $\overline{Y}$ and sample mean estimator of $\overline{X}$ based on a systematic sample of size $n$, are respectively given by

$$\overline{y}^* = \frac{n_1 \overline{y}_{n1} + n_2 \overline{y}_{h_2}}{n}$$

and $\overline{x} = \frac{1}{n}\sum_{j=1}^{n} x_{ij}$

where $\overline{y}_{n1}$ and $\overline{y}_{h_2}$ are respectively the means based on $n_1$ respondent units and $h_2$ non-respondent units. Obviously, $\overline{y}^*$ and $\overline{x}$ are unbiased estimators of $\overline{Y}$ and $\overline{X}$ respectively. The respective variances of $\overline{y}^*$ and $\overline{x}$ are expressed as

$$V\left(\overline{y}^*\right) = \frac{N-1}{nN}\{1 + (n-1)\rho_Y\}S_Y^2 + \frac{L-1}{n}W_2 S_{Y2}^2 \tag{2.1}$$

and

$$V\left(\overline{x}\right) = \frac{N-1}{nN}\{1 + (n-1)\rho_X\}S_X^2 \tag{2.2}$$

where $\rho_Y$ and $\rho_X$ are the correlation coefficients between a pair of units within the systematic sample for the study and auxiliary variables respectively. $S_Y^2$ and $S_X^2$ are respectively the mean squares of the entire group for study and auxiliary variables. $S_{Y2}^2$ be the population mean square of non-response group under study variable and $W_2$ is the non-response rate in the population.

Assuming population mean $\overline{X}$ of auxiliary variable is known, the usual ratio, product and dual to ratio estimators based on a systematic sample under non-response are respectively given by

$$\overline{y}_R^* = \frac{\overline{y}^*}{\overline{x}} \overline{X}, \tag{2.3}$$

$$\overline{y}_P^* = \frac{\overline{y}^* \overline{x}}{\overline{X}} \tag{2.4}$$

and

$$\overline{y}_D^* = \overline{y}^* \frac{\left(N\overline{X} - n\overline{x}\right)}{\left(N-n\right)\overline{X}}. \tag{2.5}$$

Obviously, all the above estimators $\overline{y}_R^*$, $\overline{y}_P^*$ and $\overline{y}_D^*$ are biased. To derive the biases and mean square errors (MSE) of the estimators $\overline{y}_R^*$, $\overline{y}_P^*$ and $\overline{y}_D^*$ under large sample approximation, let

$\overline{y}^* = \overline{Y}\left(1+e_0\right)$

$\overline{x} = \overline{X}\left(1+e_1\right)$

such that $E\left(e_0\right) = E\left(e_1\right) = 0$,

$$E\left(e_0^2\right) = \frac{V\left(\overline{y}^*\right)}{\overline{Y}^2} = \frac{N-1}{nN}\left\{1+\left(n-1\right)\rho_Y\right\}C_Y^2 + \frac{L-1}{n}W_2\frac{S_{Y2}^2}{\overline{Y}^2}, \tag{2.6}$$

$$E\left(e_1^2\right) = \frac{V\left(\overline{x}\right)}{\overline{X}^2} = \frac{N-1}{nN}\left\{1+\left(n-1\right)\rho_X\right\}C_X^2 \tag{2.7}$$

and

$$E(e_0 e_1) = \frac{Cov\left(\bar{y}^*, \bar{x}\right)}{\bar{Y}\bar{X}} = \frac{N-1}{nN}\{1+(n-1)\rho_Y\}^{1/2}\{1+(n-1)\rho_X\}^{1/2}\rho C_Y C_X \tag{2.8}$$

where $C_Y$ and $C_X$ are the coefficients of variation of study and auxiliary variables respectively in the population (for proof see Singh and Singh(1998) and Singh (2003, pg. no. 138) ).

The biases and MSE's of the estimators $\bar{y}_R^*$ , $\bar{y}_P^*$ and $\bar{y}_D^*$ up to the first order of approximation using (2.6-2.8), are respectively given by

$$B\left(\bar{y}_R^*\right) = \frac{N-1}{nN}\bar{Y}\{1+(n-1)\rho_X\}\left(1-K\rho^*\right)C_X^2, \tag{2.9}$$

$$MSE\left(\bar{y}_R^*\right) = \frac{N-1}{nN}\bar{Y}^2\{1+(n-1)\rho_X\}\left[\rho^{*2}C_Y^2+\left(1-2K\rho^*\right)C_X^2\right]+\frac{L-1}{n}W_2 S_{Y2}^2, \tag{2.10}$$

$$B\left(\bar{y}_P^*\right) = \frac{N-1}{nN}\bar{Y}\{1+(n-1)\rho_X\}K\rho^* C_X^2, \tag{2.11}$$

$$MSE\left(\bar{y}_P^*\right) = \frac{N-1}{nN}\bar{Y}^2\{1+(n-1)\rho_X\}\left[\rho^{*2}C_Y^2+\left(1+2K\rho^*\right)C_X^2\right]+\frac{L-1}{n}W_2 S_{Y2}^2, \tag{2.12}$$

$$B\left(\bar{y}_D^*\right) = \frac{\mathrm{N}-1}{\mathrm{nN}}\bar{\mathrm{Y}}\{1+(\mathrm{n}-1)\rho_{\mathrm{X}}\}\left[-\rho^*\mathrm{K}\right]\mathrm{C}_{\mathrm{X}}^2, \tag{2.13}$$

$$MSE\left(\bar{y}_D^*\right) = \frac{N-1}{nN}\bar{Y}^2\{1+(n-1)\rho_X\}\left[\rho^{*2}C_Y^2+\left(\frac{f}{1-f}\right)\left\{\left(\frac{f}{1-f}\right)-2\rho^* K\right\}C_X^2\right]$$

$$+\frac{(L-1)}{n}W_2 S_{Y2}^2 \tag{2.14}$$

where,

$$\rho^* = \{1+(n-1)\rho_Y\}^{1/2} \Big/ \{1+(n-1)\rho_X\}^{1/2} \quad \text{and } K = \rho\frac{C_Y}{C_X}.$$

( for details of proof refer to Singh et al.(2012)).

The regression estimator based on a systematic sample under non-response is given by

$$\bar{y}^*_{lr} = \bar{y}^* + b(\bar{X} - \bar{x}) \tag{2.15}$$

MSE of the estimator $\bar{y}^*_{lr}$ is given by

$$\mathrm{MSE}(\bar{y}^*_{lr}) = \frac{N-1}{nN}\bar{Y}^2\{1+(n-1)\rho_X\}\left[C_Y^2 - K^2C_X^2\right]\rho^{*2} + \frac{(L-1)}{n}W_2S_{Y2}^2 \tag{2.16}$$

## 3. Adapted Family of Estimators

Adapting the estimator proposed by Singh and Shukla (1987), a family of factor-type estimators of population mean in systematic sampling under non-response is written as

$$T^*_\alpha = \bar{y}^*\left[\frac{(A+C)\bar{X} + fB\bar{x}}{(A+fB)\bar{X} + C\bar{x}}\right]. \tag{3.1}$$

The constants A, B, C, and $f$ are same as defined in (1.1).

It can easily be seen that the proposed family generates the non-response versions of some well known estimators of population mean in systematic sampling on putting different choices of $\alpha$. For example, if we take $\alpha$ = 1, 2, 3 and 4, the resulting estimators will be ratio, product, dual to ratio and sample mean estimators of population mean in systematic sampling under non-response respectively.

### 3.1 Properties of $T^*_\alpha$

Obviously, the proposed family is biased for the population mean $\bar{Y}$. In order to find the bias and MSE of $T^*_\alpha$**,** we use large sample approximations. Expressing the equation (3.1) in terms of $e_i$'s $(i = 0,1)$ we have

$$T^*_\alpha = \frac{\bar{Y}(1+e_0)(1+De_1)^{-1}\left[(A+C)+fB(1+e_1)\right]}{A+fB+C} \tag{3.2}$$

where $D = \dfrac{C}{A+fB+C}$.

Since $|D| < 1$ and $|e_i| < 1$, neglecting the terms of $e_i$'s $(i = 0,1)$ having power greater than two, the equation (3.2) can be written as

$$T_\alpha^* - \overline{Y} = \frac{\overline{Y}}{A + fB + C}\Big[(A+C)\{e_0 - De_1 + D^2 e_1^2 - De_0 e_1\}$$

$$+ fB\{e_0 - (D-1)e_1 + D(D-1)e_1^2 - (D-1)e_0 e_1\}\Big]. \qquad (3.3)$$

Taking expectation of both sides of the equation (3.3), we get

$$E\left[T_\alpha^* - \overline{Y}\right] = \frac{\overline{Y}(C - fB)}{A + fB + C}\left[\frac{C}{A + fB + C}E\left(e_1^2\right) - E\left(e_0 e_1\right)\right].$$

Let $\phi_1(\alpha) = \dfrac{fB}{A + fB + C}$ and $\phi_2(\alpha) = \dfrac{C}{A + fB + C}$ then

$$\phi(\alpha) = \phi_2(\alpha) - \phi_1(\alpha) = \frac{C - fB}{A + fB + C}.$$

Thus, we have

$$E\left[T_\alpha^* - \overline{Y}\right] = \overline{Y}\phi(\alpha)\left[\phi_2(\alpha)E\left(e_1^2\right) - E\left(e_0 e_1\right)\right]. \qquad (3.4)$$

Putting the values of $E\left(e_1^2\right)$ and $E\left(e_0 e_1\right)$ from equations (2.7) and (2.8) into the equation (3.4), we get the bias of $T_\alpha^*$ as

$$B\left(T_\alpha^*\right) = \phi(\alpha)\frac{N-1}{nN}\overline{Y}\{1 + (n-1)\rho_X\}\left[\phi_2(\alpha) - \rho^* K\right]C_X^2. \qquad (3.5)$$

Squaring both the sides of the equation (3.3) and then taking expectation, we get

$$E\left[T_\alpha^* - \overline{Y}\right]^2 = \overline{Y}^2\left[E\left(e_0^2\right) + \phi^2(\alpha)E\left(e_1^2\right) - 2\phi(\alpha)E\left(e_0 e_1\right)\right]. \qquad (3.6)$$

Substituting the values of $E\left(e_0^2\right)$, $E\left(e_1^2\right)$ and $E\left(e_0 e_1\right)$ from the respective equations (2.6), (2.7) and (2.8) into the equation (3.6), we get the MSE of $T_\alpha^*$ as

$$MSE\left(T_\alpha^*\right) = \frac{N-1}{nN}\overline{Y}^2\{1 + (n-1)\rho_X\}\left[\rho^{*2} C_Y^2 + \left\{\phi^2(\alpha) - 2\phi(\alpha)\rho^* K\right\}C_X^2\right]$$

$$+ \frac{(L-1)}{n}W_2 S_{Y2}^2. \qquad (3.7)$$

### 3.2 Optimum Choice of $\alpha$

In order to obtain the optimum choice of $\alpha$, we differentiate the equation (3.7) with respect to $\alpha$ and equating the derivative to zero, we get the normal equation as

$$\frac{N-1}{nN}\overline{Y}^2\{1+(n-1)\rho_X\}[2\phi(\alpha)\phi'(\alpha)-2\phi'(\alpha)\rho^* K]C_X^2 = 0 \tag{3.8}$$

where $\phi'(\alpha)$ is the first derivative of $\phi(\alpha)$ with respect to $\alpha$.

Now from equation (3.8), we get

$$\phi(\alpha) = \rho^* K \tag{3.9}$$

which is the cubic equation in $\alpha$. Thus $\alpha$ has three real roots for which the MSE of proposed family would attain its minimum.

Putting the value of $\phi(\alpha)$ from equation (3.9) into equation (3.7), we get

$$MSE(T_\alpha^*)_{\min} = \frac{N-1}{nN}\overline{Y}^2\{1+(n-1)\rho_X\}[C_Y^2 - K^2 C_X^2]\rho^{*2} + \frac{(L-1)}{n}W_2 S_{Y2}^2 \tag{3.10}$$

which is the MSE of the usual regression estimator of population mean in systematic sampling under non-response.

## 4. Empirical Study

In the support of theoretical results, we have considered the data given in Murthy (1967, p. 131-132). These data are related to the length and timber volume for ten blocks of the blacks mountain experimental forest. The value of intraclass correlation coefficients $\rho_X$ and $\rho_Y$ have been given approximately equal by Murthy (1967, p. 149) and Kushwaha and Singh (1989) for the systematic sample of size 16 by enumerating all possible systematic samples after arranging the data in ascending order of strip length. The particulars of the population are given below:

$N = 176$, $n = 16$, $\overline{Y} = 282.6136$, $\overline{X} = 6.9943$,

$S_Y^2 = 24114.6700$, $S_X^2 = 8.7600$, $\rho = 0.8710$,

$$S_{Y2}^2 = \frac{3}{4} S_Y^2 = 18086.0025.$$

Table 1 depicts the MSE's and variance of the estimators of proposed family with respect to non-response rate ($W_2$).

**Table 1**: MSE and Variance of the Estimators for $L = 2$.

| $\alpha$ | $W_2$ | | | |
|---|---|---|---|---|
| | 0.1 | 0.2 | 0.3 | 0.4 |
| $1\ (=\bar{y}_R^*)$ | 371.37 | 484.41 | 597.45 | 710.48 |
| $2\ (=\bar{y}_P^*)$ | 1908.81 | 2021.85 | 2134.89 | 2247.93 |
| $3(=\bar{y}_D^*)$ | 1063.22 | 1176.26 | 1289.30 | 1402.33 |
| $4(=\bar{y}^*)$ | 1140.69 | 1253.13 | 1366.17 | 1479.205 |
| $\alpha_{opt}(=(T_\alpha^*)_{\min})$ | 270.67 | 383.71 | 496.75 | 609.78 |

## 5. Conclusion

In this paper, we have adapted Singh and Shukla (1987) estimator in systematic sampling in the presence of non-response using an auxiliary variable and obtained the optimum estimator of the proposed family. It is observed that the proposed family can generate the non-response versions of a number of estimators of population mean in systematic sampling on different choice of $\alpha$. From Table 1, we observe that the proposed family under optimum condition has minimum MSE, which is equal to the MSE of the regression estimator (most of the class of estimators in sampling literature under optimum condition attains MSE equal to the MSE of the regression estimator). It is also seen that the MSE or variance of the estimators increases with increase in non response rate in the population.

# Some Improved Estimators of Population Mean Using Information on Two Auxiliary Attributes

**[1]Hemant Verma, †[1]Rajesh Singh and [2]Florentin Smarandache**

**[1]Department of Statistics, Banaras Hindu University,Varanasi-221005, India**

**[2]Chair of Department of Mathematics, University of New Mexico, Gallup, USA**

*†Corresponding author*

## Abstract

In this paper, we have studied the problem of estimating the finite population mean when information on two auxiliary attributes are available. Some improved estimators in simple random sampling without replacement have been suggested and their properties are studied. The expressions of mean squared error's (MSE's) up to the first order of approximation are derived. An empirical study is carried out to judge the best estimator out of the suggested estimators.

## Introduction

The role of auxiliary information in survey sampling is to increase the precision of estimators when study variable is highly correlated with auxiliary variable. But when we talk about qualitative phenomena of any object then we use auxiliary attributes instead of auxiliary variable. For example, if we talk about height of a person then sex will be a good auxiliary attribute and similarly if we talk about particular breed of cow then in this case milk produced by them will be good auxiliary variable.

Most of the times, we see that instead of one auxiliary variable we have information on two auxiliary variables e.g.; to estimate the hourly wages we can use the information on marital status and region of residence (see Gujrati and Sangeetha (2007), page-311).

In this paper, we assume that both auxiliary attributes have significant point bi-serial correlation with the study variable and there is significant phi-correlation (see Yule (1912)) between the auxiliary attributes.

Consider a sample of size n drawn by simple random sampling without replacement (SRSWOR) from a population of size N. let $y_j$, $\phi_{ij}$ (i=1,2) denote the observations on variable y and $\phi_i$ (i=1,2) respectively for the $j^{th}$ unit (i=1,2,3,……N) . We note that $\phi_{ij}$=1, if $j^{th}$ unit possesses attribute $\phi_{ij}$=0 otherwise . Let $A_i = \sum_{j=1}^{N} \phi_{ij}$, $a_i = \sum_{j=1}^{n} \phi_{ij}$ ; i=1,2 denotes the total number of units in the population and sample respectively, possessing attribute $\phi$ . Similarly, let $P_i = \frac{A_i}{N}$ and $p_i = \frac{a_i}{n}$;(i=1,2 ) denotes the proportion of units in the population and sample respectively possessing attribute $\phi_i$ (i=1,2).

In order to have an estimate of the study variable y, assuming the knowledge of the population proportion P, Naik and Gupta (1996) and Singh et al. (2007) respectively proposed following estimators:

$$t_1 = \bar{y}\left(\frac{P_1}{p_1}\right) \tag{1.1}$$

$$t_2 = \bar{y}\left(\frac{p_2}{P_2}\right) \tag{1.2}$$

$$t_3 = \bar{y}\exp\left(\frac{P_1 - p_1}{P_1 + p_1}\right) \tag{1.3}$$

$$t_4 = \bar{y}\exp\left(\frac{p_2 - P_2}{p_2 + P_2}\right) \tag{1.4}$$

The bias and MSE expression's of the estimator's $t_i$ (i=1, 2, 3, 4) up to the first order of approximation are, respectively, given by

$$B(t_1) = \bar{Y} f_1 C_{p_1}^2 \left[1 - K_{pb_1}\right] \tag{1.5}$$

$$B(t_2) = \bar{Y} f_1 K_{pb_2} C_{p_2}^2 \tag{1.6}$$

$$B(t_3) = \bar{Y} f_1 \frac{C_{p_1}^2}{2}\left[\frac{1}{4} - K_{pb_1}\right] \qquad (1.7)$$

$$B(t_4) = \bar{Y} f_1 \frac{C_{p_2}^2}{2}\left[\frac{1}{4} + K_{pb_2}\right] \qquad (1.8)$$

$$MSE(t_1) = \bar{Y}^2 f_1\left[C_y^2 + C_{p_1}^2\left(1 - 2K_{pb_1}\right)\right] \qquad (1.9)$$

$$MSE(t_2) = \bar{Y}^2 f_1\left[C_y^2 + C_{p_1}^2\left(1 + 2K_{pb_2}\right)\right] \qquad (1.10)$$

$$MSE(t_3) = \bar{Y}^2 f_1\left[C_y^2 + C_{p_1}^2\left(\frac{1}{4} - K_{pb_1}\right)\right] \qquad (1.11)$$

$$MSE(t_4) = \bar{Y}^2 f_1\left[C_y^2 + C_{p_2}^2\left(\frac{1}{4} + K_{pb_2}\right)\right] \qquad (1.12)$$

where, $f_1 = \frac{1}{n} - \frac{1}{N}$ , $S_{\phi_j}^2 = \frac{1}{N-1}\sum_{i=1}^{N}\left(\phi_{ji} - P_j\right)^2$, $S_{y\phi_j} = \frac{1}{N-1}\sum_{i=1}^{N}\left(y_i - \bar{Y}\right)\left(\phi_{ji} - P_j\right)$,

$\rho_{pb_j} = \frac{S_{y\phi_j}}{S_y S_{\phi_j}}$, $C_y = \frac{S_y}{\bar{Y}}$, $C_{pj} = \frac{S_{\phi_j}}{P_j}$; $(j = 1,2)$, $K_{pb_1} = \rho_{pb_1}\frac{C_y}{C_{p_1}}$, $K_{pb_2} = \rho_{pb_2}\frac{C_y}{C_{p_2}}$.

$s_{\phi_1\phi_2} = \frac{1}{n-1}\sum_{i=1}^{n}\left(\phi_{1i} - p_1\right)\left(\phi_{2i} - p_2\right)$ and $\rho_\phi = \frac{s_{\phi_1\phi_2}}{s_{\phi_1}s_{\phi_2}}$ be the sample phi-covariance and phi-correlation between $\phi_1$ and $\phi_2$ respectively, corresponding to the population phi-covariance and phi-correlation $S_{\phi_1\phi_2} = \frac{1}{N-1}\sum_{i=1}^{N}\left(\phi_{1i} - P_1\right)\left(\phi_{2i} - P_2\right)$ and $\rho_\phi = \frac{S_{\phi_1\phi_2}}{S_{\phi_1}S_{\phi_2}}$.

In this paper we have proposed some improved estimators of population mean using information on two auxiliary attributes in simple random sampling without replacement. A comparative study is also carried out to compare the optimum estimators with respect to usual mean estimator with the help of numerical data.

## 2. Proposed Estimators

Following Olkin (1958), we propose an estimator $t_1$ as

$$t_5 = \bar{y}\left[w_1 \frac{P_1}{p_1} + w_2 \frac{P_2}{p_2}\right] \tag{2.1}$$

where $w_1$ and $w_2$ are constants, such that $w_1 + w_2 = 1$.

Consider another estimator $t_6$ as

$$t_6 = \left[K_{61}\bar{y} + K_{62}(P_1 - p_1)\right]\exp\left[\frac{P_2 - p_2}{P_2 + p_2}\right] \tag{2.2}$$

where $K_{61}$ and $K_{62}$ are constants.

Following Shaoo et al. (1993), we propose another estimator $t_7$ as

$$t_7 = \bar{y} + K_{71}(P_1 - p_1) + K_{72}(P_2 - p_2) \tag{2.3}$$

where $K_{71}$ and $K_{72}$ are constants.

**Bias and MSE of estimators** $t_5, t_6$ and $t_7$:

To obtain the bias and MSE expressions of the estimators $t_i (i = 5,6,7)$ to the first degree of approximation, we define

$$e_0 = \frac{\bar{y} - \bar{Y}}{\bar{Y}}, \quad e_1 = \frac{p_1 - P_1}{P_1}, \quad e_2 = \frac{p_2 - P_2}{P_2}$$

such that, $E(e_i) = 0$; $i = 0, 1, 2.$

Also,

$$E(e_0^2) = f_1 C_y^2, \; E(e_1^2) = f_1 C_{p_1}^2, \; E(e_2^2) = f_1 C_{p_2}^2,$$

$$E(e_0 e_1) = f_1 K_{pb_1} C_{p_1}^2, \; E(e_0 e_2) = f_1 K_{pb_2} C_{p_2}^2, \quad E(e_1 e_2) = f_1 K_{\phi} C_{p_2}^2,$$

$$K_{pb_1} = \rho_{pb_1}\frac{C_y}{C_{p_1}}, \; K_{pb_2} = \rho_{pb_2}\frac{C_y}{C_{p_2}}, \; K_{\phi} = \rho_{\phi}\frac{C_{p_1}}{C_{p_2}}$$

Expressing (2.1) in terms of e's we have,

$$t_5 = \bar{Y}(1 + e_0)\left[w_1 \frac{P_1}{P_1(1 + e_1)} + w_2 \frac{P_2}{P_2(1 + e_2)}\right]$$

$$t_5 = \overline{Y}(1+e_0)\left[w_1(1+e_1)^{-1} + w_2(1+e_2)^{-1}\right] \tag{3.1}$$

Expanding the right hand side of (3.1) and retaining terms up to second degrees of e's, we have,

$$t_5 = \overline{Y}\left[1 + e_0 - w_1 e_1 - w_2 e_2 + w_1 e_1^2 + w_2 e_2^2 - w_1 e_0 e_1 - w_2 e_0 e_2\right] \tag{3.2}$$

Taking expectations of both sides of (3.1) and then subtracting $\overline{Y}$ from both sides, we get the bias of estimator $t_5$ upto the first order of approximation as

$$\text{Bias}(t_5) = \overline{Y} f_1\left[w_1 C_{p_1}^2\left(1 - K_{pb_1}\right) + w_2 C_{p_2}^2\left(1 - K_{pb_2}\right)\right] \tag{3.3}$$

From (3.2), we have,

$$\left(t_5 - \overline{Y}\right) \cong \overline{Y}\left[e_0 - w_1 e_1 - w_2 e_2\right] \tag{3.4}$$

Squaring both sides of (3.4) and then taking expectations, we get the MSE of $t_5$ up to the first order of approximation as

$$\text{MSE}(t_5) = \overline{Y}^2 f_1\left[C_y^2 + w_1^2 C_{p_1}^2 + w_2^2 C_{p_2}^2 - 2w_1 K_{pb_1} C_{p_1}^2 - 2w_2 K_{pb_2} C_{pb_2}^2 + 2w_1 w_2 K_{\phi} C_{p_2}^2\right] \tag{3.5}$$

Minimization of (3.5) with respect to $w_1$ and $w_2$, we get the optimum values of $w_1$ and $w_2$ , as

$$w_{1(opt)} = \frac{K_{pb_1} C_{p_1}^2 - K_{\phi} C_{p_2}^2}{C_{p_1}^2 - K_{\phi} C_{p_2}^2} = w_1^*(\text{say})$$

$$\begin{aligned} w_{2(opt)} &= 1 - w_{1(opt)} \\ &= 1 - \frac{K_{pb_1} C_{p_1}^2 - K_{\phi} C_{p_2}^2}{C_{p_1}^2 - K_{\phi} C_{p_2}^2} \\ &= \frac{C_{p_1}^2\left[1 - K_{pb_1}\right]}{C_{p_1}^2 - K_{\phi} C_{p_2}^2} = w_2^*(\text{say}) \end{aligned}$$

Similarly, we get the bias and MSE expressions of estimator $t_6$ and $t_7$ respectively, as

$$\text{Bias}(t_6) = K_{61}\overline{Y}\left[1 + f_1 C_{p_2}^2\left(\frac{3}{8} - \frac{1}{2}K_{pb_2}\right)\right] + \frac{1}{2}K_{22} P_1 f_1 K_{\phi} C_{p_2}^2 \tag{3.6}$$

$$\text{Bias}(t_7) = 0 \tag{3.7}$$

And

$$MSE(t_6) = K_{61}^2 \overline{Y}^2 A_1 + K_{62}^2 P_1^2 A_2 - 2K_{61}K_{62}P_1\overline{Y}A_3 + (1-2K_{61})\overline{Y}^2 \quad (3.8)$$

where $A_1 = 1 + f_1\left(C_y^2 + C_{p_2}^2\left(\frac{1}{4} - K_{pb_2}\right)\right)$

$A_2 = f_1 C_{p_1}^2$

$A_3 = f_1\left(k_{pb_1}C_{p_1}^2 - \frac{1}{2}K_\phi C_p^2\right)$

And the optimum values of $K_{61}$ and $K_{62}$ are respectively, given as

$$K_{61(opt)} = \frac{A_2}{A_1A_2 - A_3^2} = K_{61}^*(\text{say})$$

$$K_{62(opt)} = \frac{\overline{Y}A_3}{P_1(A_1A_2 - A_3^2)} = K_{62}^*(\text{say})$$

$$MSE(t_7) = \overline{Y}^2 f_1 C_y^2 + K_{71}^2 P_1^2 f_1 C_{p_1}^2 + K_{72}^2 P_2^2 f_1 C_{p_2}^2 - 2K_{71}P_1\overline{Y}f_1K_{pb_1}C_{p_1}^2 - 2K_{72}P_2\overline{Y}f_1K_{pb_2}C_{p_2}^2 + 2K_{71}K_{72}P_1P_2f_1K_\phi C_{p_2}^2 \quad (3.9)$$

And the optimum values of $K_{71}$ and $K_{72}$ are respectively, given as

$$K_{71(opt)} = \frac{\overline{Y}}{P_1}\left(\frac{K_{pb_1}C_{p_1}^2 - K_{pb_2}K_\phi C_{p_2}^2}{C_{p_1}^2 - K_\phi^2 C_{p_2}^2}\right) = K_{71}^*(\text{say})$$

$$K_{72(opt)} = \frac{\overline{Y}}{P_2}\left(\frac{K_{pb_2}C_{p_1}^2 - K_{pb_1}K_\phi C_{p_1}^2}{C_{p_1}^2 - K_\phi^2 C_{p_2}^2}\right) = K_{72}^*(\text{say})$$

## 4. Empirical Study

Data: (Source: Government of Pakistan (2004))

The population consists rice cultivation areas in 73 districts of Pakistan. The variables are defined as:

Y= rice production (in 000' tonnes, with one tonne = 0.984 ton) during 2003,

$P_1$ = production of farms where rice production is more than 20 tonnes during the year 2002, and

$P_2$ = proportion of farms with rice cultivation area more than 20 hectares during the year 2003.

For this data, we have

N=73, $\overline{Y}$=61.3, $P_1$=0.4247, $P_2$=0.3425, $S_y^2$=12371.4, $S_{\phi_1}^2$=0.225490, $S_{\phi_2}^2$=0.228311,

$\rho_{pb_1}$=0.621, $\rho_{pb_2}$=0.673, $\rho_\phi$=0.889.

The percent relative efficiency (PRE's) of the estimators $t_i$ (i=1,2,…7) with respect to unusual unbiasedestimator $\bar{y}$ have been computed and given in Table 4.1.

**Table 4.1 : PRE of the estimators with respect to $\bar{y}$**

| Estimator | PRE |
|---|---|
| $\bar{y}$ | 100.00 |
| $t_1$ | 162.7652 |
| $t_2$ | 48.7874 |
| $t_3$ | 131.5899 |
| $t_4$ | 60.2812 |
| $t_5$ | 165.8780 |
| $t_6$ | 197.7008 |
| $t_7$ | 183.2372 |

**Conclusion**

In this paper we have proposed some improved estimators of population mean using information on two auxiliary attributes in simple random sampling without replacement. From the Table 4.1 we observe that the estimator $t_6$ is the best followed by the estimator $t_7$ .

# Study of Some Improved Ratio Type Estimators Under Second Order Approximation

**[1]Prayas Sharma, †[1]Rajesh Singh and [2]Florentin Smarandache**

**[1]Department of Statistics, Banaras Hindu University,Varanasi-221005, India**

**[2]Chair of Department of Mathematics, University of New Mexico, Gallup, USA**

***†Corresponding author***

**Abstract**

Chakrabarty(1979), Khoshnevisan et al. (2007), Sahai and Ray (1980), Ismail et al. (2011) and Solanki et al. (2012) proposed estimators for estimating population mean $\overline{Y}$. Up to the first order of approximation and under optimum conditions, the minimum mean squared error (MSE) of all the above estimators is equal to the MSE of the regression estimator. In this paper, we have tried to found out the second order biases and mean square errors of these estimators using information on auxiliary variable based on simple random sampling. Finally, we have compared the performance of these estimators with some numerical illustration.

## 1. Introduction

Let U= ($U_1$ ,$U_2$ , $U_3$, …..,$U_i$, ….$U_N$ ) denotes a finite population of distinct and identifiable units. For estimating the population mean $\overline{Y}$ of a study variable Y, let us consider X be the auxiliary variable that are correlated with study variable Y, taking the corresponding values of the units. Let a sample of size n be drawn from this population using simple random sampling without replacement (SRSWOR) and $y_i$ , $x_i$ (i=1,2,…..n ) are the values of the study variable and auxiliary variable respectively for the i-th unit of the sample.

In sampling theory the use of suitable auxiliary information results in considerable reduction in MSE of the ratio estimators. Many authors suggested estimators using some known population parameters of an auxiliary variable. Upadhyaya and Singh (1999), Singh and Tailor (2003), Kadilar and Cingi (2006), Khoshnevisan et al. (2007), Singh et al. (2007), Singh et al. (2008) and Singh and Kumar (2011) suggested estimators in simple random sampling. Most of the authors discussed the properties of estimators along with their first order bias and MSE. Hossain et al. (2006) studied some estimators in second order approximation. In this study we have studied properties of some estimators under second order of approximation.

## 2. Some Estimators in Simple Random Sampling

For estimating the population mean $\overline{Y}$ of Y, Chakrabarty (1979) proposed ratio type estimator -

$$t_1 = (1-\alpha)\overline{y} + \alpha\overline{y}\frac{\overline{X}}{\overline{x}} \tag{2.1}$$

where $\overline{y} = \frac{1}{n}\sum_{i=1}^{n} y_i$ and $\overline{x} = \frac{1}{n}\sum_{i=1}^{n} x_i$.

Khoshnevisan et al. (2007) ratio type estimator is given by

$$t_2 = \overline{y}\left[\frac{\overline{X}}{\beta\overline{x} + (1-\beta)\overline{X}}\right]^g \tag{2.2}$$

where $\beta$ and g are constants.

Sahai and Ray (1980) proposed an estimator $t_3$ as

$$t_3 = \overline{y}\left[2 - \left\{\frac{\overline{x}}{\overline{X}}\right\}^W\right] \tag{2.3}$$

Ismail et al. (2011) proposed and estimator t4 for estimating the population mean $\overline{Y}$ of Y as

$$t_4 = \bar{y}\left[\frac{\bar{x} + a(\bar{X} - \bar{x})}{\bar{x} + b(\bar{X} - \bar{x})}\right]^p \tag{2.4}$$

where p, a and b are constant.

Also, for estimating the population mean $\bar{Y}$ of Y, Solanki et al. (2012) proposed an estimator $t_5$ as

$$t_5 = \bar{y}\left[2 - \left\{\left(\frac{\bar{x}}{\bar{X}}\right)^{\lambda} \exp\frac{\delta(\bar{x} - \bar{X})}{(\bar{x} + \bar{X})}\right\}\right] \tag{2.5}$$

where $\lambda$ and $\delta$ are constants, suitably chosen by minimizing mean square error of the estimator $t_5$ .

**3. Notations used**

Let us define, $e_0 = \frac{\bar{y} - \bar{Y}}{\bar{Y}}$ and $e_1 = \frac{\bar{x} - \bar{X}}{\bar{X}}$, then $E(e_0) = E(e_1) = 0$.

For obtaining the bias and MSE the following lemmas will be used:

**Lemma 3.1**

(i) $V(e_0) = E\{(e_0)^2\} = \frac{N-n}{N-1}\frac{1}{n}C_{02} = L_1C_{02}$

(ii) $V(e_1) = E\{(e_1)^2\} = \frac{N-n}{N-1}\frac{1}{n}C_{20} = L_1C_{20}$

(iii) $COV(e_0, e_1) = E\{(e_0e_1)\} = \frac{N-n}{N-1}\frac{1}{n}C_{11} = L_1C_{11}$

**Lemma 3.2**

(iv) $E\{(e_1{}^2e_0)\} = \frac{(N-n)}{(N-1)}\frac{(N-2n)}{(N-2)}\frac{1}{n^2}C_{21} = L_2C_{21}$

(v) $E\{(e_1{}^3)\} = \frac{(N-n)}{(N-1)}\frac{(N-2n)}{(N-2)}\frac{1}{n^2}C_{30} = L_2C_{30}$

**Lemma 3.3**

(vi) $E(e_1{}^3e_0\ ) = L_3C_{31} + 3L_4C_{20}C_{11}$

(vii) $E\{(e_1{}^4)\} = \frac{(N-n))(N^2 + N - 6nN + 6n^2)}{(N-1)(N-2)(N-3)}\frac{1}{n^3}C_{30} = L_3C_{40} + 3L_4C_{20}{}^2$

(viii) $E(e_1{}^2e_0{}^2) = L_3C_{40} + 3L_4C_{20}$

Where $L_3 = \frac{(N-n))(N^2 + N - 6nN + 6n^2)}{(N-1)(N-2)(N-3)}\frac{1}{n^3}$ , $L_4 = \frac{N(N-n))(N-n-1)(n-1)}{(N-1)(N-2)(N-3)}\frac{1}{n^3}$

and $Cpq = \frac{(Xi - \overline{X})^p}{\overline{X}^p} \frac{(Y_i - \overline{Y})^q}{\overline{Y}^q}$.

Proof of these lemma's are straight forward by using SRSWOR (see Sukhatme and Sukhatme (1970)).

## 4. First Order Biases and Mean Squared Errors

The expression for the biases of the estimators $t_1$, $t_2$, $t_3$, $t_4$ and $t_5$ are respectively given by

$$\text{Bias}(t_1) = \overline{Y}\left[\frac{1}{2}\alpha L_1 C_{20} - \alpha L_1 C_{11}\right] \qquad (4.1)$$

$$\text{Bias}(t_2) = \overline{Y}\left[\frac{g(g+1)}{2} L_1 C_{20} - g\beta L_1 C_{11}\right] \qquad (4.2)$$

$$\text{Bias}(t_3) = \overline{Y}\left[-\frac{w(w-1)}{2} L_1 C_{20} - w L_1 C_{11}\right] \qquad (4.3)$$

$$\text{Bias}(t_4) = \overline{Y}\left[bDL_1 C_{20} - DL_1 C_{11} + \frac{D(b-a)(p-1)}{2} L_1 C_{20}\right] \qquad (4.4)$$

$$\text{Bias}(t_5) = \overline{Y}\left[-\frac{K(K-1)}{2} L_1 C_{20} - K L_1 C_{11}\right] \qquad (4.5)$$

where, D =p (b-a) and $k = \frac{(\delta + 2\lambda)}{2}$.

Expression for the MSE's of the estimators $t_1$, $t_2$ $t_3$ $t_4$ and $t_5$ are, respectively given by

$$\text{MSE}(t_1) = \overline{Y}^2\left[L_1 C_{02} + \alpha^2 L_1 C_{20} - 2\alpha L_1 C_{11}\right] \qquad (4.6)$$

$$\text{MSE}(t_2) = \overline{Y}^2\left[L_1 C_{02} + g^2\beta^2 L_1 C_{20} - 2g\beta L_1 C_{11}\right] \qquad (4.7)$$

$$\text{MSE}(t_3) = \overline{Y}^2\left[L_1 C_{02} + w^2 L_1 C_{20} - 2w L_1 C_{11}\right] \qquad (4.8)$$

$$\text{MSE}(t_4) = \overline{Y}^2\left[L_1 C_{02} + DL_1 C_{20} - 2DL_1 C_{11}\right] \qquad (4.9)$$

$$\text{MSE}(t_5) = \overline{Y}^2\left[L_1 C_{02} + k^2 L_1 C_{20} - 2k L_1 C_{11}\right] \qquad (4.10)$$

## 5. Second Order Biases and Mean Squared Errors

Expressing estimator $t_i$'s (i=1,2,3,4) in terms of e's (i=0,1), we get

$$t_1 = \overline{Y}(1+e_0)\left\{(1-\alpha) + \alpha(1+e_1)^{-1}\right\}$$

Or

$$t_1 - \overline{Y} = \overline{Y}\left\{e_0 + \frac{e_1}{2} + \frac{\alpha}{2}e_1{}^2 - \alpha e_0 e_1 + \alpha e_0 e_1{}^2 - \frac{\alpha}{6}e^3 - \frac{\alpha}{6}e_0 e_1{}^3 + \frac{\alpha}{24}e^4\right\} \qquad (5.1)$$

Taking expectations, we get the bias of the estimator $t_1$ up to the second order of approximation as

$$\text{Bias}_2(t_1) = \overline{Y}\left[\frac{\alpha}{2}L_1C_{20} - \alpha L_1C_{11} - \frac{\alpha}{6}L_2C_{30} + \alpha L_2C_{21} - \frac{\alpha}{6}\left(L_3C_{31} + 3L_4C_{20}C_{11}\right)\right.$$
$$\left. + \frac{\alpha}{24}(L_3C_{40} + 3L_4C_{20}{}^2)\right] \qquad (5.2)$$

Similarly, we get the biases of the estimator's $t_2$, $t_3$, $t_4$ and $t_5$ up to second order of approximation, respectively as

$$\text{Bias}_2(t_2) = \overline{Y}\left[\frac{g(g+1)}{2}\beta^2L_1C_{20} - g\beta L_1C_{11} - \frac{g(g+1)}{2}\beta^2L_2C_{21} - \frac{g(g+1)(g+2)}{6}\beta^3L_2C_{30}\right.$$
$$- \frac{g(g+1)(g+2)}{6}\beta^3\left(L_3C_{31} + 3L_4C_{20}C_{11}\right)$$
$$\left. + \frac{g(g+1)(g+2)(g+3)}{24}\beta^4(L_3C_{40} + 3L_4C_{20}{}^2)\right] \qquad (5.3)$$

$$\text{Bias}_2(t_3) = \overline{Y}\left[-\frac{w(w-1)}{2}L_1C_{20} - wL_1C_{11} - \frac{w(w-1)}{2}L_2C_{21} - \frac{w(w-1)(w-2)}{6}L_2C_{30}\right.$$
$$- \frac{w(w-1)(w-2)}{6}\left(L_3C_{31} + 3L_4C_{20}C_{11}\right)$$
$$\left. - \frac{w(w-1)(w-2)(w-3)}{24}(L_3C_{40} + 3L_4C_{20}{}^2)\right] \qquad (5.4)$$

$$\text{Bias}_2(t_4) = \overline{Y}\left[\frac{(D_1bD)}{2}L_1C_{20} - DL_1C_{11} + \frac{(bD + D_1)}{}L_2C_{21} - \frac{(b^2D + 2bD_1 + D_2)}{}L_2C_{30}\right.$$
$$- (b^2D + 2bD_1)\left(L_3C_{31} + 3L_4C_{20}C_{11}\right)$$
$$\left. + \frac{(b^3D + 3b^2D_1 + 3bD_2 + D_3)}{}(L_3C_{40} + 3L_4C_{20}{}^2)\right] \qquad (5.5)$$

where, $D_1 = D\frac{(b-a)(p-1)}{2}$ $\quad D_2 = D_1\frac{(b-a)(p-2)}{3}$.

$$\text{Bias}_2(t_5) = \overline{Y}\Big[-\frac{k(k-1)}{2}L_1C_{20} - kL_1C_{11} - \frac{k(k-1)}{2}L_2C_{21} - ML_2C_{30} - M(L_3C_{31} + 3L_4C_{20}C_{11})$$

$$-N(L_3C_{40} + 3L_4C_{20}{}^2)\Big] \qquad (5.6)$$

Where, $M=\frac{1}{2}\left\{\frac{(\delta^3 - 6\delta^2)}{24} + \frac{(\alpha(\delta^2 - 2\delta)}{4} + \frac{\lambda(\lambda-1)}{2}\delta + \frac{\lambda(\lambda-1)(\lambda-2)}{3}\right\}$,

$k = \frac{(\delta + 2\lambda)}{2}$,

$N=\frac{1}{8}\left\{\frac{(\delta^4 - 12\delta^3 + 12\delta^2)}{48} + \frac{(\alpha(\delta^3 - 6\delta)}{6} + \frac{\lambda(\lambda-1)}{2}(\delta^2 - 2\delta) + \frac{\lambda(\lambda-1)(\lambda-2)(\lambda-3)}{3}\right\}$.

The MSE's of the estimators $t_1$, $t_2$, $t_3$, $t_4$ and $t_5$ up to the second order of approximation are, respectively given by

$$\text{MSE}_2(t_1) = \overline{Y}^2\Big[L_1C_{02} + \alpha^2L_1C_{20} - 2\alpha L_1C_{11} - \alpha^2L_2C_{30} + (2\alpha^2 + \alpha)L_2C_{21}$$
$$-2\alpha^2(L_3C_{31} + 3L_4C_{20}C_{11})$$
$$+\alpha(\alpha+1)(L_3C_{22} + 3L_4(C_{20}C_{02} + C_{11}^2)) + \frac{5}{24}\alpha^2(L_3C_{40} + 3L_4C_{20}{}^2)\Big] \qquad (5.7)$$

$$\text{MSE}_2(t_2) = \overline{Y}^2\Big[L_1C_{02} + g^2\beta^2L_1C_{20} - 2\beta gL_1C_{11} - \beta^3g^2(g+1)L_2C_{30} + g(3g+1)\beta^2L_2C_{21}$$
$$-2\beta gL_2C_{12} - \left\{\frac{7g^3 + 9g^2 + 2g}{3}\right\}\beta^3(L_3C_{31} + 3L_4C_{20}C_{11})$$
$$+g(2g+1)\beta^2(L_3C_{22} + 3L_4(C_{20}C_{02} + C_{11}^2))$$
$$+\left\{\frac{2g^3 + 9g^2 + 10g + 3}{6}\right\}\beta^4(L_3C_{40} + 3L_4C_{20}{}^2\Big] \qquad (5.8)$$

$$MSE_2(t_3)=\overline{Y}^2\Big[L_1C_{02}+w^2L_1C_{20}-2wL_1C_{11}-w^2(w-1)L_2C_{30}+w(w+1)L_2C_{21}-2wL_2C_{12}$$
$$+\left\{\frac{5w^3-3w^2-2w}{3}\right\}\left(L_3C_{31}+3L_4C_{20}C_{11}\right)$$
$$+w\left(L_3C_{22}+3L_4(C_{20}C_{02}+C_{11}^2)\right)+\left\{\frac{7w^4-18w^3+11w^2}{24}\right\}(L_3C_{40}+3L_4C_{20}{}^2)\Big] \tag{5.9}$$

$$MSE_2(t_4)=\overline{Y}^2\Big[L_1C_{02}+D^2L_1C_{20}-2DL_1C_{11}-4DD_1L_2C_{30}+(2bD+2D_1+2D^2)L_2C_{21}-2DL_2C_{12}$$
$$+\left\{2D^2+2b^2D+2DD_1+4bD_1+4bD^2\right\}\left(L_3C_{31}+3L_4C_{20}C_{11}\right)$$
$$+\left\{D^2+2D_1+2bD\right\}\left(L_3C_{22}+3L_4(C_{20}C_{02}+C_{11}^2)\right)$$
$$+\left\{3b^2D^2+D_1^2+2DD_2+12bDD_1\right\}(L_3C_{40}+3L_4C_{20}^2)\Big] \tag{5.10}$$

$$MSE_2(t_5)=\overline{Y}^2\Big[L_1C_{02}+k^2L_1C_{20}-2kL_1C_{11}+kL_2C_{21}-2kL_2C_{12}+k^2(k-1)L_2C_{30}$$
$$+2k^2(k-1)\left(L_3C_{31}+3L_4C_{20}C_{11}\right)+k\left(L_3C_{22}+3L_4(C_{20}C_{02}+C_{11}^2)\right)$$
$$+\frac{(k^2-k)^2}{4}(L_3C_{40}+3L_4C_{20}{}^2)\Big] \tag{5.11}$$

## 6. Numerical Illustration

For a natural population data, we have calculated the biases and the mean square error's of the estimator's and compare these biases and MSE's of the estimator's under first and second order of approximations.

**<u>Data Set</u>**

The data for the empirical analysis are taken from 1981, Utter Pradesh District Census Handbook, Aligar. The population consist of 340 villages under koil police station, with Y=Number of agricultural labour in 1981 and X=Area of the villages (in acre) in 1981. The following values are obtained

$\overline{Y} = 73.76765,\ \overline{X} = 2419.04,\ \ N = 340, n = 70, n' = 120$, n=70, $C_{02}$=0.7614, $C_{11}$=0.2667, $C_{03}$=2.6942, $C_{12}$=0.0747, $C_{12}$=0.1589, $C_{30}$=0.7877, $C_{13}$=0.1321, $C_{31}$=0.8851, $C_{04}$=17.4275 $C_{22}$=0.8424, $C_{40}$=1.3051

**Table 6.1: Biases and MSE's of the estimators**

| Estimator | Bias | | MSE | |
|---|---|---|---|---|
| | **First order** | **Second order** | **First order** | **Second order** |
| $t_1$ | 0.0044915 | 0.004424 | 39.217225 | 39.45222 |
| $t_2$ | 0 | -0.00036 | 39.217225 (for g=1) | 39.33552 (for g=1) |
| $t_3$ | -0.04922 | -0.04935 | 39.217225 | 39.29102 |
| $t_4$ | 0.2809243 | -0.60428 | 39.217225 | 39.44855 |
| $t_5$ | -0.027679 | -0.04911 | 39.217225 | 39.27187 |

In the Table 6.1 the biases and MSE's of the estimators $t_1$, $t_2$, $t_3$, $t_4$ and $t_5$ are written under first order and second order of approximations. For all the estimators $t_1$, $t_2$, $t_3$, $t_4$ and $t_5$, it was observed that the value of the biases decreased and the value of the MSE's increased for second order approximation. MSE's of the estimators up to the first order of approximation under optimum conditions are same. From Table 6.1 we observe that under second order of approximation the estimator $t_5$ is best followed by $t_3$,and $t_2$ among the estimators considered here for the given data set.

## 7. Estimators under stratified random sampling

In survey sampling, it is well established that the use of auxiliary information results in substantial gain in efficiency over the estimators which do not use such information. However, in planning surveys, the stratified sampling has often proved needful in improving the precision of estimates over simple random sampling. Assume that the population U consist of L strata as $U=U_1, U_2,\ldots,U_L$ . Here the size of the stratum $U_h$ is $N_h$, and the size of simple random sample in stratum $U_h$ is $n_h$, where h=1, 2,---,L.

The Chakrabarty(1979) ratio- type estimator under stratified random sampling is given by

$$t_1' = (1-\alpha)\bar{y}_{st} + \alpha\bar{y}_{st}\frac{\bar{X}}{\bar{x}_{st}} \tag{7.1}$$

where,

$$\bar{y}_h = \frac{1}{n_h}\sum_{i=1}^{n_h} y_{hi}, \qquad \bar{x}_h = \frac{1}{n_h}\sum_{i=1}^{n_h} x_{hi},$$

$$\bar{y}_{st} = \sum_{h=1}^{L} w_h\bar{y}_h, \quad \bar{x}_{st} = \sum_{h=1}^{L} w_h\bar{x}_h, \text{ and } \quad \bar{X} = \sum_{h=1}^{L} w_h\bar{X}_h.$$

Khoshnevisan et al. (2007) ratio- type estimator under stratified random sampling is given by

$$t_2' = \bar{y}_{st}\left[\frac{\bar{X}}{\beta\bar{x}_{st} + (1-\beta)\bar{X}}\right]^g \tag{7.2}$$

where g is a constant, for g=1 , $t_2'$ is same as conventional ratio estimator whereas for g = -1, it becomes conventional product type estimator.

Sahai and Ray (1980) estimator $t_3$ under stratified random sampling is given by

$$t_3' = \bar{y}_{st}\left[2 - \left\{\frac{\bar{x}_{st}}{\bar{X}}\right\}^W\right] \tag{7.3}$$

Ismail et al. (2011) estimator under stratified random sampling $t_4'$ is given by

$$t'_4 = \bar{y}_{st}\left[\frac{\bar{x} + a(\bar{X} - \bar{x}_{st})}{\bar{x} + b(\bar{X} - \bar{x}_{st})}\right]^p \tag{7.4}$$

Solanki et al. (2012) estimator under stratified random sampling is given as

$$t'_5 = \bar{y}_{st}\left[2 - \left\{\left(\frac{\bar{x}_{st}}{\bar{X}}\right)^{\lambda} \exp\frac{\delta(\bar{x}_{st} - \bar{X})}{(\bar{x}_{st} + \bar{X})}\right\}\right] \tag{7.5}$$

where $\lambda$ and $\delta$ are the constants, suitably chosen by minimizing MSE of the estimator $t'_5$.

## 8. Notations used under stratified random sampling

Let us define, $e_0 = \frac{\bar{y}_{st} - \bar{y}}{\bar{y}}$ and $e_1 = \frac{\bar{x}_{st} - \bar{x}}{\bar{x}}$, then $E(e_0) = E(e_1) = 0$.

To obtain the bias and MSE of the proposed estimators, we use the following notations in the rest of the article:

$$\bar{y}_{st} = \sum_{h=1}^{L} w_h \bar{y}_h = \bar{Y}(1 + e_0),$$

$$\bar{x}_{st} = \sum_{h=1}^{L} w_h \bar{x}_h = \bar{X}(1 + e_1),$$

such that,

$E(e_0) = E(e_1) = E(e_2) = 0,$ and

$$V_{rs} = \sum_{h=1}^{L} W_h^{r+s} E\left[(\bar{x}_h - \bar{X}_h)^r (\bar{y}_h - \bar{Y}_h)^s\right]$$

Also

$$E(e_0^2) = \frac{\sum_{h=1}^{L} w_h^2 \gamma_h S_{yh}^2}{\bar{Y}^2} = V_{20}$$

$$E(e_1^2) = \frac{\sum_{h=1}^{L} w_h^2 \gamma_h S_{xyh}^2}{\bar{X}^2} = V_{02}$$

$$E(e_0e_1)=\frac{\sum_{h=1}^{L} w_h^2\gamma_h S_{xyh}^2}{\overline{X}\overline{Y}}=V_{11}$$

Where

$$S_{yh}^2=\frac{\sum_{i=1}^{N_h}(\overline{y}_h-\overline{Y}_h)^2}{N_h-1}, \quad S_{xh}^2=\frac{\sum_{i=1}^{N_h}(\overline{x}_h-\overline{X}_h)^2}{N_h-1}, \quad S_{xyh}=\frac{\sum_{i=1}^{N_h}(\overline{x}_h-\overline{X}_h)(\overline{y}_h-\overline{Y}_h)}{N_h-1}$$

$$\gamma_h=\frac{1-f_h}{n_h}, \qquad f_h=\frac{n_h}{N_h}, \qquad w_h=\frac{N_h}{n_h}.$$

Some additional notations for second order approximation,

$$V_{rs}=\sum_{h=1}^{L} W_h^{r+s}\frac{1}{\overline{Y}^r\overline{X}^s}E\left[\left(\overline{y}_h-\overline{Y}_h\right)^s\left(\overline{x}_h-\overline{X}_h\right)^r\right]$$

Where,

$$C_{rs(h)}=\frac{1}{N_h}\sum_{i=1}^{N_h}\left[\left(\overline{y}_h-\overline{Y}_h\right)^s\left(\overline{x}_h-\overline{X}_h\right)^r\right]$$

$$V_{12}=\sum_{h=1}^{L} W_h^3\frac{k_{1(h)}C_{12(h)}}{\overline{Y}\overline{X}^2} \qquad V_{21}=\sum_{h=1}^{L} W_h^3\frac{k_{1(h)}C_{21(h)}}{\overline{Y}^2\overline{X}} \qquad V_{30}=\sum_{h=1}^{L} W_h^3\frac{k_{1(h)}C_{30(h)}}{\overline{Y}^3}$$

$$V_{03}=\sum_{h=1}^{L} W_h^3\frac{k_{1(h)}C_{03(h)}}{\overline{X}^3} \qquad V_{13}=\sum_{h=1}^{L} W_h^4\frac{k_{2(h)}C_{13(h)}+3k_{3(h)}C_{01(h)}C_{02(h)}}{\overline{Y}\overline{X}^3}$$

$$V_{04}=\sum_{h=1}^{L} W_h^4\frac{k_{2(h)}C_{04(h)}+3k_{3(h)}C_{02(h)}^2}{\overline{X}^4}$$

$$V_{22}=\sum_{h=1}^{L} W_h^4\frac{k_{2(h)}C_{22(h)}+k_{3(h)}\left(C_{01(h)}C_{02(h)}+2C_{11(h)}^2\right)}{\overline{Y}^2\overline{X}^2}$$

Where

$$k_{1(h)} = \frac{(N_h - n_h)(N_h - 2n_h)}{n^2(N_h - 1)(N_h - 2)}$$

$$k_{2(h)} = \frac{(N_h - n_h)(N_h + 1)N_h - 6n_h(N_h - n_h)}{n^3(N_h - 1)(N_h - 2)(N_h - 3)}$$

$$k_{3(h)} = \frac{(N_h - n_h)N_h(N_h - n_h - 1)(n_h - 1)}{n^3(N_h - 1)(N_h - 2)(N_h - 3)}$$

## 9. First Order Biases and Mean Squared Errors

The biases of the estimators $t'_1, t'_2, t'_3, t'_4$ and $t'_5$ are respectively given by

$$\text{Bias}(t'_1) = \overline{Y}\left[\frac{1}{2}\alpha L_1 V_{02} - \alpha V_{11}\right] \tag{9.1}$$

$$\text{Bias}(t'_2) = \overline{Y}\left[\frac{g(g+1)}{2}\beta^2 V_{02} - g\beta V_{11}\right] \tag{9.2}$$

$$\text{Bias}(t'_3) = \overline{Y}\left[-\frac{w(w-1)}{2}V_{02} - wV_{11}\right] \tag{9.3}$$

$$\text{Bias}(t'_4) = \overline{Y}\left[bDV_{02} - DV_{11} + \frac{D(b-a)(p-1)}{2}V_{02}\right] \tag{9.4}$$

$$\text{Bias}(t'_5) = \overline{Y}\left[-\frac{K(K-1)}{2}V_{02} - KV_{11}\right] \tag{9.5}$$

Where, $D = p(b-a)$ and $k = \frac{(\delta + 2\lambda)}{2}$.

The MSE's of the estimators $t'_1, t'_2, t'_3, t'_4$ and $t'_5$ are respectively given by

$$\text{MSE}(t'_1) = \overline{Y}^2\left[V_{20} + \alpha^2 V_{02} - 2\alpha V_{11}\right] \tag{9.6}$$

$$\text{MSE}(t'_2) = \overline{Y}^2\left[V_{20} + g^2\beta^2 V_{02} - 2g\beta V_{11}\right] \tag{9.7}$$

$$\text{MSE}(t'_3) = \overline{Y}^2\left[V_{20} + w^2 V_{02} - 2wV_{11}\right] \tag{9.8}$$

$$\text{MSE}(t'_4) = \overline{Y}^2\left[V_{20} + DV_{02} - 2DV_{11}\right] \tag{9.9}$$

$$\text{MSE}(t'_5) = \overline{Y}^2\left[V_{20} + k^2 V_{02} - 2kV_{11}\right] \tag{9.10}$$

## 10. Second Order Biases and Mean Squared Errors

Expressing estimator $t_i$'s (i=1,2,3,4) in terms of e's (i=0,1), we get

$$t_1' = \overline{Y}(1+e_0)\{(1-\alpha)+\alpha(1+e_1)^{-1}\}$$

Or

$$t_1' - \overline{Y} = \overline{Y}\left\{e_0 + \frac{e_1}{2} + \frac{\alpha}{2}e_1^{\ 2} - \alpha e_0 e_1 + \alpha e_0 e_1^{\ 2} - \frac{\alpha}{6}e^3 - \frac{\alpha}{6}e_0 e_1^{\ 3} + \frac{\alpha}{24}e^4\right\} \qquad (10.1)$$

Taking expectations, we get the bias of the estimator $t_1'$ up to the second order of approximation as

$$Bias_2(t_1') = \overline{Y}\left[\frac{\alpha}{2}V_{02} - \alpha V_{11} - \frac{\alpha}{6}V_{03} + \alpha V_{12} - \frac{\alpha}{6}V_{13} + \frac{\alpha}{24}V_{04}\right] \qquad (10.2)$$

Similarly we get the Biases of the estimator's $t_2', t_3', t_4'$ and $t_5'$ up to second order of approximation, respectively as

$$Bias_2(t_2') = \overline{Y}\left[\frac{g(g+1)}{2}\beta^2 V_{02} - g\beta V_{11} - \frac{g(g+1)}{2}\beta^2 V_{12} - \frac{g(g+1)(g+2)}{6}\beta^3 V_{30}\right.$$
$$\left. - \frac{g(g+1)(g+2)}{6}\beta^3 V_{31} + \frac{g(g+1)(g+2)(g+3)}{24}\beta^4 V_{04}\right] \qquad (10.3)$$

$$Bias_2(t_3') = \overline{Y}\left[-\frac{w(w-1)}{2}V_{02} - wV_{11} - \frac{w(w-1)}{2}V_{12} - \frac{w(w-1)(w-2)}{6}V_{03}\right.$$
$$\left. - \frac{w(w-1)(w-2)}{6}V_{31} - \frac{w(w-1)(w-2)(w-3)}{24}V_{04}\right] \qquad (10.4)$$

$$Bias_2(t_4') = \overline{Y}\left[\frac{(D_1 bD)}{2}V_{02} - DV_{11} + (bD + D_1)V_{12} - (b^2 D + 2bD_1 + D_2)V_{03}\right.$$
$$\left. - (b^2 D + 2bD_1)V_{31} + (b^3 D + 3b^2 D_1 + 3bD_2 + D_3)V_{04}\right] \qquad (10.5)$$

Where, $D = p(b-a)$ $\qquad D_1 = D\frac{(b-a)(p-1)}{2}$ $\qquad D_2 = D_1\frac{(b-a)(p-2)}{3}$

$$Bias_2(t_5') = \overline{Y}\left[\left[-\frac{k(k-1)}{2}V_{02} - kV_{11} - \frac{k(k-1)}{2}V_{12} - MV_{03} - MV_{31} - NV_{04}\right]\right. \qquad (10.6)$$

Where, $M=\frac{1}{2}\left\{\frac{\left(\delta^3-6\delta^2\right)}{24}+\frac{\left(\alpha(\delta^2-2\delta\right)}{4}+\frac{\lambda(\lambda-1)}{2}\delta+\frac{\lambda(\lambda-1)(\lambda-2)}{3}\right\}$, $k=\frac{(\delta+2\lambda)}{2}$

$$N=\frac{1}{8}\left\{\frac{\left(\delta^4-12\delta^3+12\delta^2\right)}{48}+\frac{\left(\alpha(\delta^3-6\delta\right)}{6}+\frac{\lambda(\lambda-1)}{2}(\delta^2-2\delta)+\frac{\lambda(\lambda-1)(\lambda-2)(\lambda-3)}{3}\right\}$$

Following are the MSE of the estimators $t'_1, t'_2, t'_3, t'_4$ and $t'_5$ up to second order of approximation

$$MSE_2(t'_1)=\overline{Y}^2\Big[V_{20}+\alpha^2V_{02}-2\alpha V_{11}-\alpha^2V_{03}+(2\alpha^2+\alpha)V_{12}$$
$$-2\alpha^2V_{31}+\alpha(\alpha+1)V_{22}+\frac{5}{24}\alpha^2V_{04}\Big] \tag{10.7}$$

$$MSE_2(t'_2)=\overline{Y}^2\Big[V_{20}+g^2\beta^2V_{02}-2\beta gV_{11}-\beta^3g^2(g+1)V_{03}+g(3g+1)\beta^2V_{12}$$
$$-2\beta gV_{21}-\left\{\frac{7g^3+9g^2+2g}{3}\right\}\beta^3V_{31}+g(2g+1)\beta^2V_{22}$$
$$+\left\{\frac{2g^3+9g^2+10g+3}{6}\right\}\beta^4V_{04}\Big] \tag{10.8}$$

$$MSE_2(t'_3)=\overline{Y}^2\Big[V_{20}+w^2V_{02}-2wV_{11}-w^2(w-1)V_{03}+w(w+1)V_{12}-2wV_{21}$$
$$+\left\{\frac{5w^3-3w^2-2w}{3}\right\}V_{31}+wV_{22}+\left\{\frac{7w^4-18w^3+11w^2}{24}\right\}V_{04}\Big] \tag{10.9}$$

$$MSE_2(t'_4)=\overline{Y}^2\Big[V_{20}+D^2V_{02}-2DV_{11}-4DD_1V_{03}+(2bD+2D_1+2D^2)V_{12}-2DV_{21}$$
$$+\left\{2D^2+2b^2D+2DD_1+4bD_1+4bD^2\right\}V_{31}+\left\{D^2+2D_1+2bD\right\}V_{22}$$
$$+\left\{3b^2D^2+D_1^2+2DD_2+12bDD_1\right\}V_{04}\Big] \tag{10.10}$$

$$MSE_2(t'_5)=\overline{Y}^2\Big[V_{20}+k^2V_{02}-2kV_{11}+kV_{12}-2kV_{21}+k^2(k-1)V_{03}$$
$$+2k^2(k-1)V_{31}+kV_{22}+\frac{(k^2-k)^2}{4}V_{04}\Big] \tag{10.11}$$

### 11. Numerical Illustration

For the natural population data, we shall calculate the bias and the mean square error of the estimator and compare Bias and MSE for the first and second order of approximation.

**Data Set-1**

To illustrate the performance of above estimators, we have considered the natural Data given in *Singh and Chaudhary (1986, p.162)*.

The data were collected in a pilot survey for estimating the extent of cultivation and production of fresh fruits in three districts of Uttar- Pradesh in the year 1976-1977.

**Table 11.1: Biases and MSE's of the estimators**

| Estimator | Bias | | MSE | |
|---|---|---|---|---|
| | First order | Second order | First order | second order |
| $t_1'$ | -10.82707903 | -13.65734654 | 1299.110219 | 1372.906438 |
| $t_2'$ | -10.82707903 | 6.543275811 | 1299.110219 | 1367.548263 |
| $t_3'$ | -27.05776113 | -27.0653128 | 1299.110219 | 1417.2785 |
| $t_4'$ | 11.69553975 | -41.84516913 | 1299.110219 | 2605.736045 |
| $t_5'$ | -22.38574093 | -14.95110301 | 1299.110219 | 2440.644397 |

From Table 11.1 we observe that the MSE's of the estimators $t'_1, t'_2, t'_3, t'_4$ and $t'_5$ are same up to the first order of approximation but the biases are different. The MSE of the estimator $t'_2$ is minimum under second order of approximation followed by the estimator $t'_1$ and other estimators.

## Conclusion

In this study we have considered some estimators whose MSE's are same up to the first order of approximation. We have extended the study to second order of approximation to search for best estimator in the sense of minimum variance. The properties of the estimators are studied under SRSWOR and stratified random sampling. We have observed from Table 6.1 and Table 11.1 that the behavior of the estimators changes dramatically when we consider the terms up to second order of approximation.

# IMPROVEMENT IN ESTIMATING THE POPULATION MEAN USING DUAL TO RATIO-CUM-PRODUCT ESTIMATOR IN SIMPLE RANDOM SAMPLING

**[1]Olufadi Yunusa, [2]†Rajesh Singh and [3]Florentin Smarandache**

**[1]Department of Statistics and Mathematical Sciences**

**Kwara State University, P.M.B 1530, Malete, Nigeria**

**[2]Department of Statistics, Banaras Hindu University, Varanasi(U.P.), India**

**[3]Chair of Department of Mathematics, University of New Mexico, Gallup, USA**

***†Corresponding author***

## ABSTRACT

In this paper, we propose a new estimator for estimating the finite population mean using two auxiliary variables. The expressions for the bias and mean square error of the suggested estimator have been obtained to the first degree of approximation and some estimators are shown to be a particular member of this estimator. Furthermore, comparison of the suggested estimator with the usual unbiased estimator and other estimators considered in this paper is carried out. In addition, an empirical study with two natural data from literature is used to expound the performance of the proposed estimator with respect to others.

## 1. INTRODUCTION

It is well known that the use of auxiliary information in sample survey design results in efficient estimate of population parameters (e.g. mean) under some realistic conditions. This information may be used at the design stage (leading, for instance, to stratification, systematic or probability proportional to size sampling designs), at the estimation stage or at both stages. The literature on survey sampling describes a great variety of techniques for using auxiliary information by means of ratio, product and regression methods. Ratio and product type estimators take advantage of the correlation between the auxiliary variable, $x$ and the study variable, $y$. For example, when information is available on the auxiliary variable that is positively (high) correlated with the study variable, the ratio method of estimation is a suitable estimator to estimate the population mean and when the correlation is negative the product method of estimation as envisaged by Robson (1957) and Murthy (1964) is appropriate.

Quite often information on many auxiliary variables is available in the survey which can be utilized to increase the precision of the estimate. In this situation, Olkin (1958) was the first author to deal with the problem of estimating the mean of a survey variable when auxiliary variables are made available. He suggested the use of information on more than one supplementary characteristic, positively correlated with the study variable, considering a linear combination of ratio estimators based on each auxiliary variable separately. The coefficients of the linear combination were determined so as to minimize the variance of the estimator. Analogously to Olkin, Singh (1967) gave a multivariate expression of Murthy's (1964) product estimator, while Raj (1965) suggested a method for using multi-auxiliary variables through a linear combination of single difference estimators. More recently, Abu-Dayyeh et al. (2003), Kadilar and Cingi (2004, 2005), Perri (2004, 2005), Dianna and Perri (2007), Malik and Singh (2012) among others have suggested estimators for $\bar{Y}$ using information on several auxiliary variables.

Motivated by Srivenkataramana (1980), Bandyopadhyay (1980) and Singh et al. (2005) and with the aim of providing a more efficient estimator; we propose, in this paper, a new estimator for $\bar{Y}$ when two auxiliary variables are available.

## 2. BACKGROUND TO THE SUGGESTED ESTIMATOR

Consider a finite population $P = (P_1, P_2, ..., P_N)$ of $N$ units. Let a sample $s$ of size $n$ be drawn from this population by simple random sampling without replacements (SRSWOR). Let $y_i$ and $(x_i, z_i)$ represents the value of a response variable $y$ and two auxiliary variables $(x, z)$ are available. The units of this finite population are identifiable in the sense that they are uniquely labeled from 1 to $N$ and the label on each unit is known. Further, suppose in a survey problem, we are interested in estimating the population mean $\bar{Y}$ of $y$, assuming that the population means $(\bar{X}, \bar{Z})$ of $(x, z)$ are known. The traditional ratio and product estimators for $\bar{Y}$ are given as

$$\bar{y}_R = \bar{y}\frac{\bar{X}}{\bar{x}} \text{ and } \bar{y}_P = \bar{y}\frac{\bar{z}}{\bar{Z}}$$

respectively, where $\bar{y} = \frac{1}{n}\sum_{i=1}^{n} y_i$, $\bar{x} = \frac{1}{n}\sum_{i=1}^{n} x_i$ and $\bar{z} = \frac{1}{n}\sum_{i=1}^{n} z_i$ are the sample means of $y$, $x$ and $z$ respectively.

Singh (1969) improved the ratio and product method of estimation given above and suggested the "ratio-cum-product" estimator for $\bar{Y}$ as $\bar{y}_S = \bar{y}\frac{\bar{X}}{\bar{x}}\frac{\bar{z}}{\bar{Z}}$

In literature, it has been shown by various authors; see for example, Reddy (1974) and Srivenkataramana (1978) that the bias and the mean square error of the ratio estimator $\bar{y}_R$, can be reduced with the application of transformation on the auxiliary variable $x$. Thus, authors like, Srivenkataramana (1980), Bandyopadhyay (1980) Tracy et al. (1996), Singh et al. (1998), Singh et al. (2005), Singh et al. (2007), Bartkus and Plikusas (2009) and Singh et al. (2011) have improved on the ratio, product and ratio-cum-product method of estimation using the transformation on the auxiliary information. We give below the transformations employed by these authors:

$$x_i^* = (1+g)\bar{X} - gx_i \text{ and } z_i^* = (1+g)\bar{Z} - gz_i \text{, for } i = 1,2,...,N, \quad (1)$$

where $g = \dfrac{n}{N-n}$.

Then clearly, $\bar{x}^* = (1+g)\bar{X} - g\bar{x}$ and $\bar{z}^* = (1+g)\bar{Z} - g\bar{z}$ are also unbiased estimate of $\bar{X}$ and $\bar{Z}$ respectively and $Corr(\bar{y}, \bar{x}^*) = -\rho_{yx}$ and $Corr(\bar{y}, \bar{z}^*) = -\rho_{yz}$. It is to be noted that by using the transformation above, the construction of the estimators for $\bar{Y}$ requires the knowledge of unknown parameters, which restrict the applicability of these estimators. To overcome this restriction, in practice, information on these parameters can be obtained approximately from either past experience or pilot sample survey, inexpensively.

The following estimators $\bar{y}_R^*$, $\bar{y}_P^*$ and $\bar{y}_{SE}$ are referred to as dual to ratio, product and ratio-cum-product estimators and are due to Srivenkataramana (1980), Bandyopadhyay (1980) and Singh et al. (2005) respectively. They are as presented: $\bar{y}_R^* = \bar{y}\dfrac{\bar{x}^*}{\bar{X}}$, $\bar{y}_P^* = \bar{y}\dfrac{\bar{Z}}{\bar{z}^*}$ and $\bar{y}_{SE} = \bar{y}\dfrac{\bar{x}^*}{\bar{X}}\dfrac{\bar{Z}}{\bar{z}^*}$

It is well known that the variance of the simple mean estimator $\bar{y}$, under SRSWOR design is

$$V(\bar{y}) = \lambda S_y^2$$

and to the first order of approximation, the Mean Square Errors (MSE) of $\bar{y}_R$, $\bar{y}_P$, $\bar{y}_S$, $\bar{y}_R^*$, $\bar{y}_P^*$ and $\bar{y}_{SE}$ are, respectively, given by

$$MSE(\bar{y}_R) = \lambda\left(S_y^2 + R_1^2 S_x^2 - 2R_1 S_{yx}\right)$$

$$MSE(\bar{y}_P) = \lambda\left(S_y^2 + R_2^2 S_z^2 + 2R_2 S_{yz}\right)$$

$$MSE(\bar{y}_S) = \lambda\left[S_y^2 - 2D + C\right]$$

$$MSE(\bar{y}_R^*) = \lambda\left(S_y^2 + g^2 R_1^2 S_x^2 - 2gR_1 S_{yx}\right)$$

$$MSE(\bar{y}_P^*) = \lambda\left(S_y^2 + g^2 R_2^2 S_z^2 + 2gR_2 S_{yz}\right)$$

$$MSE(\bar{y}_{SE}) = \lambda\left(S_y^2 + g^2C - 2gD\right)$$

where,

$$\lambda = \frac{1-f}{n}, \quad f = \frac{n}{N}, \quad S_y^2 = \frac{1}{N}\sum_{i=1}^{N}\left(y_i - \bar{Y}\right)^2, \quad S_{yx} = \frac{1}{N}\sum_{i=1}^{N}\left(y_i - \bar{Y}\right)\left(x_i - \bar{X}\right), \quad \rho_{yx} = \frac{S_{yx}}{S_y S_x},$$

$R_1 = \frac{\bar{Y}}{\bar{X}}$, $R_2 = \frac{\bar{Y}}{\bar{Z}}$, $C = R_1^2 S_x^2 - 2R_1R_2S_{zx} + R_2^2 S_z^2$, $D = R_1S_{yx} - R_2S_{yz}$ and $S_j^2$ for $(j = x, y, z)$ represents the variances of $x$, $y$ and $z$ respectively; while $S_{yx}$, $S_{yz}$ and $S_{zx}$ denote the covariance between $y$ and $x$, $y$ and $z$ and $z$ and $x$ respectively. Note that $\rho_{yz}$, $\rho_{zx}$, $S_x^2$, $S_z^2$, $S_{yz}$ and $S_{zx}$ are defined analogously and respective to the subscripts used.

More recently, Sharma and Tailor (2010) proposed a new ratio-cum-dual to ratio estimator of finite population mean in simple random sampling, their estimator with its MSE are respectively given as,

$$\bar{y}_{ST} = \bar{y}\left[\alpha\left(\frac{\bar{X}}{\bar{x}}\right) + (1-\alpha)\left(\frac{\bar{x}^*}{\bar{X}}\right)\right]$$

$$MSE(\bar{y}_{ST}) = \lambda S_y^2\left(1 - \rho_{yx}^2\right).$$

## 3. PROPOSED DUAL TO RATIO-CUM-PRODUCT ESTIMATOR

Using the transformation given in (1), we suggest a new estimator for $\bar{Y}$ as follows:

$$\bar{y}_{PR} = \bar{y}\left[\theta\left(\frac{\bar{x}^*}{\bar{X}}\frac{\bar{Z}}{\bar{z}^*}\right) + (1-\theta)\left(\frac{\bar{X}}{\bar{x}^*}\frac{\bar{z}^*}{\bar{Z}}\right)\right]$$

We note that when information on the auxiliary variable $z$ is not used (or variable $z$ takes the value `unity') and $\theta = 1$, the suggested estimator $\bar{y}_{PR}$ reduces to the `dual to ratio' estimator $\bar{y}_R^*$ proposed by Srivenkataramana (1980). Also, $\bar{y}_{PR}$ reduces to the `dual to product' estimator $\bar{y}_P^*$ proposed by Bandyopadhyay (1980) estimator if the information on the auxiliary variate $x$ is not used and $\theta = 0$. Furthermore, the suggested estimator reduces to the dual to ratio-cum-product estimator suggested by Singh et al. (2005) when $\theta = 1$ and information on the two auxiliary variables $x$ and $z$ are been utilized.

In order to study the properties of the suggested estimator $\bar{y}_{PR}$ (e.g. MSE), we write

$$\bar{y} = \bar{Y}(1 + k_0);\ \bar{x} = \bar{X}(1 + k_1);\ \bar{z} = \bar{Z}(1 + k_2);$$

with $E(k_0) = E(k_1) = E(k_2) = 0$

and

$$E(k_0^2) = \frac{\lambda S_y^2}{\bar{Y}^2};\quad E(k_1^2) = \frac{\lambda S_x^2}{\bar{X}^2};\quad E(k_2^2) = \frac{\lambda S_z^2}{\bar{Z}^2};\quad E(k_0 k_1) = \frac{\lambda S_{yx}}{\bar{Y}\bar{X}};\quad E(k_0 k_2) = \frac{\lambda S_{yz}}{\bar{Y}\bar{Z}};$$

$$E(k_1 k_2) = \frac{\lambda S_{zx}}{\bar{X}\bar{Z}},$$

Now expressing $\bar{y}_{PR}$ in terms of $k's$, we have

$$\bar{y}_{PR} = \bar{Y}(1 + k_0)\left[\theta(1 - gk_1)(1 - gk_2)^{-1} + (1 - \theta)(1 - gk_1)^{-1}(1 - gk_2)\right] \qquad (2)$$

We assume that $|gk_1| < 1$ and $|gk_2| < 1$ so that the right hand side of (2) is expandable. Now expanding the right hand side of (2) to the first degree of approximation, we have

$$\bar{y}_{PR} - \bar{Y} = \bar{Y}\left[k_0 + (1 - 2\alpha)g(k_1 - k_2 + k_0 k_1 - k_0 k_2) + g^2\left(k_1^2 - k_1 k_2 - \alpha\left(k_1^2 - k_2^2\right)\right)\right] \qquad (3)$$

Taking expectations on both sides of (3), we get the bias of $\bar{y}_{PR}$ to the first degree of approximation, as

$$B(\bar{y}_{PR}) = \lambda\bar{Y}\left[gDA + g^2\left(R_1^2 S_x^2 - R_1 R_2 S_{zx} - \theta\left(R_1^2 S_x^2 - R_2^2 S_z^2\right)\right)\right]$$

where $A = 1 - 2\theta$

Squaring both sides of (3) and neglecting terms of $k's$ involving power greater than two, we have

$$(\bar{y}_{PR} - \bar{Y})^2 = \bar{Y}^2\left[k_0 + Agk_1 - Agk_2\right]^2$$

$$= \bar{Y}^2\left[k_0^2 + 2Agk_0 k_1 - 2Agk_0 k_2 - 2A^2 g^2 k_1 k_2 + A^2 g^2 k_1^2 + A^2 g^2 k_2^2\right] \qquad (4)$$

Taking expectations on both sides of (4), we get the MSE of $\bar{y}_{PR}$, to the first order of approximation, as

$$MSE(\bar{y}_{PR}) = \lambda\left[S_y^2 + 2AgD + A^2 g^2 C\right] \quad (5)$$

The MSE of the proposed estimator given in (5) can be re-written in terms of coefficient of variation as

$$MSE(\bar{y}_{PR}) = \lambda\bar{Y}^2\left[C_y^2 + 2AgC_y D^* + A^2 g^2 C^*\right]$$

where $C^* = C_x^2 + C_z^2 - 2\rho_{zx}C_z C_x$ and $D^* = \rho_{yx}C_x - \rho_{yz}C_z$, $C_y = \frac{S_y}{\bar{Y}}$, $C_x = \frac{S_x}{\bar{X}}$, $C_z = \frac{S_z}{\bar{Z}}$

The MSE equation given in (5) is minimized for

$$\theta = \frac{D + Cg}{2Cg} = \theta_0 \text{ (say)}$$

We can obtain the minimum MSE of the suggested estimator $\bar{y}_{PR}$, by using the optimal equation of $\theta$ in (5) as follows: $\min.MSE(\bar{y}_{PR}) = \lambda\left[S_y^2 + F(2D + CF)\right]$

where $F = g - E$ and $E = \frac{D + Cg}{C}$

## 3. EFFICIENCY COMPARISON

In this section, the efficiency of the suggested estimator $\bar{y}_{PR}$ over the following estimator, $\bar{y}$, $\bar{y}_R$, $\bar{y}_P$, $\bar{y}_S$, $\bar{y}_R^*$, $\bar{y}_P^*$, $\bar{y}_{SE}$ and $\bar{y}_{ST}$ are investigated. We will have the conditions as follows:

(a) $MSE(\bar{y}_{PR}) - V(\bar{y}) < 0$ if $\theta > \frac{2D + gC}{2gC}$

(b) $MSE(\bar{y}_{PR}) - MSE(\bar{y}_R) < 0$ if

$$Ag(2D + AgC) < R_1\left(R_1 S_x^2 - 2S_{yx}\right) \text{ provided } S_{yx} < \frac{R_1 S_x^2}{2}$$

(c) $MSE(\bar{y}_{PR}) - MSE(\bar{y}_P) < 0$ if

$$Ag(2D+AgC)<R_2\left(R_2S_z^2+2S_{yz}\right) \text{ provided } S_{yz}<\frac{R_2S_z^2}{2}$$

(d) $MSE(\bar{y}_{PR})-MSE(\bar{y}_S)<0$ if $C<\frac{-2D}{(Ag-1)}$ provided $g<\frac{1}{A}$

(e) $MSE(\bar{y}_{PR})-MSE(\bar{y}_R^*)<0$ if

$$4\theta(\theta gC-gC-D)<2gR_1R_2S_{zx}+2R_2S_{yz}-4R_1S_{yx}-gR_2^2S_z^2$$

(f) $MSE(\bar{y}_{PR})-MSE(\bar{y}_P^*)<0$ if

$$4\theta(\theta gC-gC-D)<2gR_1R_2S_{zx}+4R_2S_{yz}-2R_1S_{yx}-gR_1^2S_x^2$$

(g) $MSE(\bar{y}_{PR})-MSE(\bar{y}_{SE})<0$ if $g<\frac{-2D}{C(A-1)}$ provided $A<1$

(h) $MSE(\bar{y}_{PR})-MSE(\bar{y}_{ST})<0$ if $Ag(2D+AgC)<-\rho^2S_y^2$

## 4. NUMERICAL ILLUSTRATION

To analyze the performance of the suggested estimator in comparison to other estimators considered in this paper, two natural population data sets from the literature are being considered. The descriptions of these populations are given below.

(1) Population I [Singh (1969, p. 377]; a detailed description can be seen in Singh (1965)

$y$ : Number of females employed

$x$ : Number of females in service

$z$ : Number of educated females

$N=61$, $n=20$, $\bar{Y}=7.46$, $\bar{X}=5.31$ , $\bar{Z}=179$, $S_y^2=28.0818$, $S_x^2=16.1761$, $S_z^2=2028.1953$, $\rho_{xy}=0.7737$, $\rho_{yz}=-0.2070$, $\rho_{zx}=-0.0033$,

(2) Population II [Source: Johnston 1972, p. 171]; A detailed description of these variables is shown in Table 1.

$y$ : Percentage of hives affected by disease

$x$ : Mean January temperature

$z$ : Date of flowering of a particular summer species (number of days from January 1)

$N=10$, $n=4$, $\bar{Y}=52$, $\bar{X}=42$ , $\bar{Z}=200$, $S_y^2=65.9776$, $S_x^2=29.9880$, $S_z^2=84$,

$\rho_{xy}=0.8$, $\rho_{yz}=-0.94$, $\rho_{zx}=-0.073$,

For these comparisons, the Percent Relative Efficiencies (PREs) of the different estimators are computed with respect to the usual unbiased estimator $\bar{y}$, using the formula

$$PRE(.,\bar{y})=\frac{V(\bar{y})}{MSE(.)}\times 100$$

and they are as presented in Table 2.

**Table 1: Description of Population II.**

| $y$ | $x$ | $z$ |
|---|---|---|
| 49 | 35 | 200 |
| 40 | 35 | 212 |
| 41 | 38 | 211 |
| 46 | 40 | 212 |
| 52 | 40 | 203 |
| 59 | 42 | 194 |
| 53 | 44 | 194 |
| 61 | 46 | 188 |
| 55 | 50 | 196 |
| 64 | 50 | 190 |

Table 2 shows clearly that the proposed dual to ratio-cum-product estimator $\bar{y}_{PR}$ has the highest PRE than other estimators; therefore, we can conclude based on the study populations that the suggested estimator is more efficient than the usual unbiased estimators, the traditional ratio and product estimator, ratio-cum-product estimator by Singh (1969),

Srivenkataramana (1980) estimator, Bandyopadhyay (1980) estimator, Singh et al. (2005) estimator and Sharma and Tailor (2010).

**Table 2: PRE of the different estimators with respect to $\bar{y}$**

| Estimators | Population I | Population II |
|---|---|---|
| $\bar{y}$ | 100 | 100 |
| $\bar{y}_R$ | 205 | 277 |
| $\bar{y}_P$ | 102 | 187 |
| $\bar{y}_S$ | 214 | 395 |
| $\bar{y}_R^*$ | 215 | 239 |
| $\bar{y}_P^*$ | 105 | 150 |
| $\bar{y}_{SE}$ | 236 | 402 |
| $\bar{y}_{ST}$ | 250 | 278 |
| $\bar{y}_{PR}$ | 279 | 457 |

## 5. CONCLUSION

We have developed a new estimator for estimating the finite population mean, which is found to be more efficient than the usual unbiased estimator, the traditional ratio and product estimators and the estimators proposed by Singh (1969), Srivenkataramana (1980), Bandyopadhyay (1980), Singh et al. (2005) and Sharma and Tailor (2010). This theoretical inference is also satisfied by the result of an application with original data. In future, we hope to extend the estimators suggested here for the development of a new estimator in stratified random sampling.

# Some Improved Estimators For Population Variance Using Two Auxiliary Variables In Double Sampling

**[1]Viplav Kumar Singh, †[1]Rajesh Singh and [2]Florentin Smarandache**

**[1]Department of Statistics, Banaras Hindu University,Varanasi-221005, India**

**[2]Chair of Department of Mathematics, University of New Mexico, Gallup, USA**

***†Corresponding author***

## Abstract

In this article we have proposed an efficient generalised class of estimator using two auxiliary variables for estimating unknown population variance $S_y^2$ of study variable y .We have also extended our problem to the case of two phase sampling. In support of theoretical results we have included an empirical study.

## 1. Introduction

Use of auxiliary information improves the precision of the estimate of parameter .Out of many ratio and product methods of estimation are good example in this context. We can use ratio method of estimation when correlation coefficient between auxiliary and study variate is positive (high), on the other hand we use product method of estimation when correlation coefficient between auxiliary and study variate is highly negative.

Variations are present everywhere in our day-to-day life. An agriculturist needs an adequate understanding of the variations in climatic factors especially from place to place (or time to time) to be able to plan on when, how and where to plant his crop. The problem of estimation of finite population variance $S_y^2$, of the study variable $y$ was discussed by Isaki (1983), Singh and Singh (2001, 2002, 2003), Singh et al. (2008), Grover (2010), and Singh et al. (2011).

Let x and z are auxiliary variates having values $(x_i, z_i)$ and y is the study variate having values $(y_i)$ respectively. Let $V_i (i = 1,2,.......N)$ is the population having N units such that y is positively correlated x and negatively correlated with z. To estimate $S_y^2$,we assume that $S_x^2$ and $S_z^2$ are known ,where

$$S_y^2 = \frac{1}{N}\sum_{i=1}^{N}(y_i - \overline{Y})^2 , \; S_x^2 = \frac{1}{N}\sum_{i=1}^{N}(x_i - \overline{X})^2 \text{ and } S_z^2 = \frac{1}{N}\sum_{i=1}^{N}(z_i - \overline{Z})^2.$$

Assume that N is large so that the finite population correction terms are ignored. A sample of size n is drawn from the population V using simple random sample without replacement.

Usual unbiased estimator of population variance $S_y^2$ is $s_y^2$, where, $s_y^2 = \frac{1}{(n-1)}\sum_{i=1}^{n}(y_i - \overline{y})^2$.

Up to the first order of approximation, variance of $s_y^2$ is given by

$$\text{var}(s_y^2) = \frac{S_y^4}{n}\partial_{400}^* \tag{1.1}$$

where, $\partial_{400}^* = \partial_{400} - 1$, $\quad \partial_{pqr} = \frac{\mu_{pqr}}{\mu_{200}^{p/2}\mu_{020}^{q/2}\mu_{002}^{r/2}}$, and

$\mu_{pqr} = \frac{1}{N}\sum_{i=1}^{N}(y_i - \overline{Y})^p (x_i - \overline{X})^q (z_i - Z)^r$ ; p,q,r being the non negative integers.

## 2. **Existing Estimators**

Let $s_y^2 = S_y^2(1 + e_0), s_x^2 = S_x^2(1 + e_1)$ and $s_z^2 = S_z^2(1 + e_2)$

where, $s_x^2 = \frac{1}{(n-1)}\sum_{i=1}^{n}(x_i - \overline{x})^2, s_z^2 = \frac{1}{(n-1)}\sum_{i=1}^{n}(z_i - \overline{z})^2$

and $\overline{x} = \frac{1}{n}\sum_{i=1}^{n}x_i, \quad \overline{z} = \frac{1}{n}\sum_{i=1}^{n}z_i$.

Also, let

$E(e_1) = E(e_2) = 0$

$$E(e_0^2)=\frac{\partial^*_{400}}{n}, E(e_1^2)=\frac{\partial^*_{040}}{n} \text{ and } E(e_2^2)=\frac{\partial^*_{004}}{n}$$

$$E(e_0e_1)=\frac{1}{n}\partial^*_{220}, E(e_1e_2)=\frac{1}{n}\partial^*_{022}, E(e_0e_2)=\frac{1}{n}\partial^*_{202}$$

Isaki (1983) suggested ratio estimator $t_1$ for estimating $S_y^2$ as-

$$t_1=s_y^2\frac{S_y^2}{s_x^2}\text{; where } s_x^2 \text{ is unbiased estimator of } S_x^2 \tag{1.2}$$

Up to the first order of approximation, mean square error of $t_1$ is given by,

$$MSE(t_1)=\frac{S_y^4}{n}\left[\partial^*_{400}+\partial^*_{040}-2\partial^*_{220}\right] \tag{1.3}$$

Singh et al. (2007) proposed the exponential ratio-type estimator $t_2$ as-

$$t_2=s_y^2\exp\left[\frac{S_x^2-s_x^2}{S_x^2+s_x^2}\right] \tag{1.4}$$

And exponential product type estimator $t_3$ as-

$$t_3=s_y^2\left[\frac{s_x^2-S_x^2}{s_x^2+S_x^2}\right] \tag{1.5}$$

Following Kadilar and Cingi (2006), Singh et al. (2011) proposed an improved estimator for estimating population variance $S_y^2$, as-

$$t_4=s_y^2\left[k_4\exp\left\{\frac{S_x^2-s_x^2}{S_x^2+s_x^2}\right\}+(1-k_4)\left\{\frac{s_x^2-S_x^2}{s_x^2+S_x^2}\right\}\right] \tag{1.6}$$

where $k_4$ is a constant.

Up to the first order of approximation mean square errors of $t_2, t_3$ and $t_4$ are respectively given by

$$MSE(t_2)=\frac{S_y^4}{n}\left[\partial^*_{400}+\frac{\partial^*_{040}}{4}-\partial^*_{220}\right] \tag{1.7}$$

$$MSE(t_3)=\frac{S_y^4}{n}\left[\partial^*_{400}+\frac{\partial^*_{004}}{4}-\partial^*_{202}\right] \tag{1.8}$$

$$MSE(t_4)=\frac{S_y^4}{n}\left[\partial^*_{400}+k_4^2\frac{\partial^*_{040}}{4}+(1-k_4)^2\frac{\partial^*_{004}}{4}-k_4\partial^*_{220}+(1-k_4)\partial^*_{202}-\frac{k_4(1-k_4)}{2}\partial^*_{022}\right] \tag{1.9}$$

where $k_4=\dfrac{(\partial^*_{004}/2)+\partial^*_{220}+\partial^*_{022}}{2(\partial^*_{040}+\partial^*_{004}+\partial^*_{022}}$.

# 3. Improved Estimator

Using Singh and Solanki (2011), we propose some improved estimators for estimating population variance $S_y^2$ as-

$$t_5=s_y^2\left[\frac{cS_x^2-Ds_x^2}{(c-d)S_x^2}\right]^p \tag{1.10}$$

$$t_6=s_y^2\left[\frac{(a+b)S_z^2}{aS_x^2+bs_x^2}\right]^q \tag{1.11}$$

$$t_7=s_y^2\left[k_7\left\{\frac{cS_x^2-Ds_x^2}{(c-d)S_x^2}\right\}^p+(1-k_7)\left\{\frac{(a+b)S_z^2}{aS_x^2+bs_x^2}\right\}^q\right] \tag{1.12}$$

where a, b, c, d are suitably choosen constants and $k_7$ is a real constant to be determined so as to minimize MSE's.

Expressing $t_5, t_6$ and $t_7$ in terms of $e_i$'s, we have

$$t_5=S_y^2\left[1-x_1pe_1+e_0-x_1pe_0e_1\right] \tag{1.13}$$

where, $x_1=\dfrac{d}{(c-d)}$.

$$t_6=S_y^2\left[1-qx_2e_2+e_0-qx_2e_0e_2\right] \tag{1.14}$$

where, $x_2=\dfrac{b}{(a+b)}$.

$$t_7=S_y^2\left[1+e_0+px_1k_7(e_1'-e_1)+(k_7-1)qx_2e_2'\right] \tag{1.15}$$

The mean squared error of estimators are obtained by subtracting $S_y^2$ from each estimator and squaring both sides and than taking expectations-

$$MSE(t_5) = \frac{S_y^4}{n}\left[\partial_{400}^* + x_1^2 p^2 \partial_{040}^* - 2x_1 p \partial_{220}^*\right] \tag{1.16}$$

Differentiating (1.16) with respect to $x_1$, we get the optimum value of $x_1$ as-

$$x_{1(opt)} = \frac{\partial_{220}^*}{p\partial_{040}^*}.$$

$$MSE(t_6) = \frac{S_y^4}{n}\left[\partial_{400}^* + x_2^2 q^2 \partial_{004}^* - 2x_1 \partial_{202}^*\right] \tag{1.17}$$

Differentiating (1.17) with respect to $x_2$, we get the optimum value of $x_2$ as –

$$x_{2(opt)} = \frac{\partial_{202}^*}{q\partial_{004}^*}.$$

$$MSE(t_7) = \frac{S_y^4}{n}\left[A + k_7^2 B + (1-k_7)C - 2k_7(1-k_7)E - 2(1-k_7)F\right] \tag{1.18}$$

Differentiating (1.18) with respect to $k_7$, we get the optimum value $k_7$ of as –

$$k_{7(opt)} = \frac{C + D - F - E}{B + C - 2E}.$$

where,

$A = \partial_{400}^*$, $B = x_1^2 p^2 \partial_{040}^*$,

$C = x_2^2 q^2 \partial_{004}^*$, $D = x_1 p \partial_{220}^*$,

$E = x_1 x_2 pq \partial_{022}^*$, $F = x_2 q \partial_{202}^*$.

## 2. Estimators In Two Phase Sampling

In certain practical situations when $S_x^2$ is not known a priori, the technique of two phase sampling or double sampling is used. Allowing SRSWOR design in each phase, the two – phase sampling scheme is as follows:

- The first phase sample $s_n^{'} (s_n^{'} \subset V)$ of a fixed size n' is drawn to measure only x and z in order to formulate the a good estimate of $S_x^2$ and $S_z^2$, respectively.

- Given $s'_n$ ,the second phase sample $s_n (s_n \subset s'_n)$ of a fixed size n is drawn to measure y only.

# Existing Estimators

Singh et al. (2007) proposed some estimators to estimate $S_y^2$ in two phase sampling, as:

$$t'_2 = s_y^2 \exp\left[\frac{s'^2_x - s^2_x}{s'^2_x + s^2_x}\right] \quad (2.1)$$

$$t'_3 = s_y^2 \exp\left[\frac{s'^2_z - s^2_z}{s'^2_z + s^2_z}\right] \quad (2.2)$$

$$t'_4 = s_y^2\left[k'_4 \exp\left\{\frac{s'^2_x - s^2_x}{s'^2_x + s^2_x}\right\} + (1 - k'_4)\exp\left\{\frac{s'^2_z - s^2_z}{s'^2_z + s^2_z}\right\}\right] \quad (2.3)$$

MSE of the estimator $t'_2, t'_3$ and $t'_4$ are respectively, given by

$$MSE(t'_2) = S_y^4\left[\frac{\partial^*_{400}}{n} + \frac{1}{4}\left(\frac{1}{n} - \frac{1}{n'}\right)\partial^*_{040} + \left(\frac{1}{n'} - \frac{1}{n}\right)\partial^*_{220}\right] \quad (2.4)$$

$$MSE(t'_3) = S_y^4\left[\frac{\partial^*_{400}}{n} + \frac{1}{4}\left(\frac{1}{n} - \frac{1}{n'}\right)\partial^*_{004} - \left(\frac{1}{n'} - \frac{1}{n}\right)\partial^*_{202}\right] \quad (2.5)$$

$$MSE(t'_4) = S_y^4\left[A' + k'^2_4 B' + (1 - k'_4)^2 C' + k'_4 D' + (1 - k'_4)E'\right] \quad (2.6)$$

And $k'_{4(opt)} = \dfrac{2C' + E' - D'}{2(B' + C')}$

Where,

$$A' = \frac{\partial^*_{400}}{n}, \quad B' = \frac{1}{4}\left(\frac{1}{n} - \frac{1}{n'}\right)\partial^*_{040} \quad C' = \frac{1}{4n'}\partial^*_{004}$$

$$D' = \left(\frac{1}{n'} - \frac{1}{n}\right)\partial^*_{202}, \quad E' = \frac{1}{n'}\partial^*_{202}$$

# Proposed estimators in two phase sampling

The estimator proposed in section 3 will take the following form in two phase sampling;

$$t'_5 = s_y^2\left[\frac{cs'^2_x - ds^2_x}{(c - d)s'^2_x}\right] \quad (2.7)$$

$$t_6^{'} = s_y^2 \left[ \frac{(a+b)S_z^2}{aS_z^2 + bs'^2_z} \right] \tag{2.8}$$

$$t_7^{'} = s_y^2 \left[ k_7^{'} \left\{ \frac{cs'^2_x - ds_x^2}{(c-d)s'^2_x} \right\} + (1-k_7^{'}) \left\{ \frac{(a+b)S_z^2}{aS_z^2 + bs'^2_z} \right\} \right] \tag{2.9}$$

Let,

$$s_y^2 = S_y^2(1+e_0)\ , s_x^2 = S_x^2(1+e_1), s'^2_x = S_x^2(1+e_1^{'})$$

$$s_z^2 = S_z^2(1+e_2)\ , s'^2_z = S_z^2(1+e_2^{'})$$

Where,

$$s'^2_x = \frac{1}{(n'-1)} \sum_{i=1}^{n'} (x_i - \bar{x}')^2 , s'^2_z = \frac{1}{(n'-1)} \sum_{i=1}^{n'} (z_i - \bar{z}')^2$$

and $\bar{x}' = \frac{1}{n'} \sum_{i=1}^{n'} x_i , \bar{z}' = \frac{1}{n'} \sum_{i=1}^{n'} z_i$

Also,

$$E(e_1^{'}) = E(e_2^{'}) = 0$$

$$E(e_1^{'2}) = \frac{\partial^*_{040}}{n'}, E(e_2^{'2}) = \frac{\partial^*_{004}}{n'}$$

$$E(e_0 e_1^{'}) = \frac{1}{n'}\partial^*_{220}\ , E(e_0 e_2^{'}) = \frac{1}{n'}\partial^*_{202}\ , E(e_1 e_1^{'}) = \frac{1}{n'}\partial^*_{040}$$

$$E(e_2 e_2^{'}) = \frac{1}{n'}\partial^*_{004}\ , E(e_1 e_2^{'}) = \frac{1}{n'}\partial^*_{022}\ , E(e_1{}' e_2^{'}) = \frac{1}{n'}\partial^*_{022}$$

Writing estimators $t_5^{'}, t_6^{'}$ and $t_7^{'}$ in terms of $e_i$'s we have ,respectively

$$t_5^{'} = S_y^2 \left[ 1 + e_0 + px_1(e_1^{'} - e_1) \right] \tag{2.10}$$

$$t_6^{'} = S_y^2 \left[ 1 + e_0 - qx_2 e_2^{'} \right] \tag{2.11}$$

$$t_7^{'} = S_y^2 \left[ 1 + px_1 k_7^{'}((e_1^{'} - e_1) + e_0 + (k_7^{'} - 1)qx_2 e_2^{'} \right] \tag{2.12}$$

Solving (1.10),(1.11) and (1.12),we get the MSE'S of the estimators $t_5^{'}$, $t_6^{'}$ and $t_7^{'}$, respectively as-

$$MSE(t_5^{'}) = S_y^4 \left[ \frac{\partial^*_{400}}{n} + p^2 x_1^2 \left( \frac{1}{n} - \frac{1}{n'} \right) \partial^*_{040} + 2px_1 \left( \frac{1}{n'} - \frac{1}{n} \right) \partial^*_{220} \right] \tag{2.13}$$

Differentiate (2.13) w.r.t. $x_1$, we get the optimum value of $x_1$ as-

$$x_{1(opt)} = \frac{\partial^*_{220}}{p\partial^*_{040}}$$

$$MSE(t'_6) = S_y^4\left[\frac{\partial^*_{400}}{n} + q^2 x_2^2 \frac{1}{n'}\partial^*_{004} + 2qx_2 \frac{1}{n'}\partial^*_{202}\right] \qquad (2.14)$$

Differentiate (2.14) with respect to $x_2$, we get the optimum value of $x_2$ as –

$$x_{2(opt)} = \frac{\partial^*_{202}}{q\partial^*_{004}}.$$

$$MSE(t'_7) = S_y^4\left[A_1 + B_1 k_7^{'2} + (k'_7 - 1)^2 C_1 + 2k'_7 D_1 + 2(k'_7 - 1)E_1\right] \qquad (2.15)$$

Differentiate (2.15) with respect to $k'_7$, we get the optimum value of $k'_7$ as –

$$k'_{7(opt)} = \frac{C_1 - D_1 - E_1}{B_1 + C_1}$$

Where,

$$A_1 = \frac{\partial^*_{400}}{n},\ B_1 = p^2 x_1^2\left(\frac{1}{n} - \frac{1}{n'}\right)\partial^*_{040},\ C_1 = px_1 \frac{q^2 x_2^2}{n'}\partial^*_{004}$$

$$D_1 = px_1\left(\frac{1}{n'} - \frac{1}{n}\right)\partial^*_{220},\ E_1 = qx_2 \frac{\partial^*_{202}}{n'}$$

# 5. Empirical Study

In support of theoretical result an empirical study is carried out. The data is taken from Murthy(1967):

$$\partial_{400} = 3.726, \partial_{040} = 2.912, \partial_{044} = 2.808$$

$$\partial_{022} = 2.73, \partial_{202} = 2.979, \partial_{220} = 3.105$$

$$c_x = 0.5938, c_y = 0.7531, c_z = 0.7205$$

$$\rho_{yz} = 0.904, \rho_{xy} = 0.98, n = 7, n' = 15$$

$$\overline{X} = 747.5882, \overline{Y} = 199.4412, \overline{Z} = 208.8824$$

❖ In Table 5.1 percent relative efficiency of various estimators of $S_y^2$ is written with respect to $s_y^2$

➢ **Table 5.1: PRE of the estimator with respect to $s_y^2$**

| Estimators | PRE |
|---|---|
| $s_y^2$ | 100 |
| $t_1$ | 636.9158 |
| $t_2$ | 248.0436 |
| $t_3$ | 52.86019 |
| $t_4$ | **699.2526** |
| $t_5$ | 667.2895 |
| $t_6$ | 486.9362 |
| $t_7$ | **699.5512** |

- In Table 5.2 percent relative efficiency of various estimators of $S_y^2$ is written with respect to $s_y^2$ in two phase sampling:

- **Table 5.2: PRE of the estimators in two phase sampling with respect to $s_y^2$**

| Estimators | PRE |
|---|---|
| $s_y^2$ | 100 |
| $t_2'$ | 142.60 |
| $t_3'$ | 66.42 |
| $t_4'$ | **460.75** |
| $t_5'$ | 182.95 |
| $t_6'$ | 158.93 |
| $t_7'$ | **568.75** |

# 6. Conclusion

In Table 5.1 and 5.2 percent relative efficiencies of various estimators are written with respect to $s_y^2$. From Table 5.1 and 5.2 we observe that the proposed estimator under optimum condition performs better than usual estimator, Isaki (1983) estimator and Singh et al. (2007 ) estimator.